\documentclass[10pt, a4paper]{article}

\usepackage{mathptmx}      
\usepackage{helvet}       
\usepackage{courier}      
\usepackage{type1cm}        
\usepackage{graphicx}        
\usepackage{multicol}        
\usepackage[bottom]{footmisc}
\usepackage{tikz}
\usetikzlibrary{calc}
\makeindex            
\usepackage{amssymb}
\usepackage{amsmath}
\usepackage{mathtools}

\def\muEVP{$\mu$\hspace{-0.3mm}EVP}

\def\muRef{{\hat{\mu}}}
\def\aRef{{\hat{a}}}
\def\eRef{{\hat{e}}}

\def\aNormRefDual#1{\left\|#1\right\|_{\muRef; X_h'}}

\def\lambdaRB#1{\lambda_{\mathrm{red},\,#1}}

\def\uRB#1{u_{\mathrm{red},\,#1}}

\def\d{d}
\DeclareMathOperator{\divergence}{div}
\newcommand{\ew}{\lambda}
\newcommand{\lm}{ {\tau} }
\newcommand{\lmtest}{ \widehat{\tau} }
\newcommand{\R}{\mathbb{R}}
\newcommand{\hatOmega}{{\widehat{\Omega}}}
\newcommand{\hatu}{{\widehat{u}}}
\newcommand{\hatv}{{\widehat{v}}}
\newcommand{\hatx}{{\widehat{\bf x}}}

\newcommand{\hatepsilon}{{\widehat{\varepsilon}}}

\newcommand{\DF}{{DF(\hatx;\mu)}}

\begin{document}
\title{Reduced basis isogeometric mortar approximations for eigenvalue problems in vibroacoustics
\thanks{We would like to gratefully acknowledge the funds provided by the “Deutsche Forschungsgemeinschaft” under the
contract/grant numbers: WO 671/11-1, WO 671/13-2 and  WO 671/15-1 (within the Priority Programme SPP 1748, "Reliable Simulation Techniques in Solid Mechanics. Development of Non-standard Discretisation Methods, Mechanical and Mathematical Analysis").} }
\author{ {\sc  Thomas Horger\thanks{Corresponding author. E-Mail: horger@ma.tum.de}, 
Barbara Wohlmuth \thanks{wohlmuth@ma.tum.de},  
 Linus Wunderlich \thanks{linus.wunderlich@ma.tum.de} }
\\[4pt] Institute for Numerical Mathematics, Technische Universit\"at M\"unchen,\\  Boltzmannstra\ss{}e 3, 85748 Garching b. M\"unchen, Germany }
\date{}
\maketitle

\begin{abstract}{We simulate the vibration of a violin bridge in a multi-query context using
reduced basis techniques. The mathematical model is based on an eigenvalue problem for the orthotropic linear elasticity equation. In addition to the nine material parameters, a geometrical  thickness
parameter is considered. This parameter enters as a 10th material parameter into the system by a mapping onto a parameter independent reference domain.  The detailed simulation is carried out by isogeometric mortar methods. Weakly coupled patch-wise tensorial structured isogeometric elements are of special interest for complex geometries with piecewise smooth but curvilinear boundaries.  To obtain locality in the detailed system, we use the saddle point approach and do not apply static condensation techniques. However within 
the reduced basis context, it is natural to eliminate the Lagrange multiplier and formulate a reduced eigenvalue problem for a symmetric positive definite matrix.  The selection of the snapshots is controlled by a multi-query greedy strategy taking into account an error indicator allowing for multiple eigenvalues.}
\end{abstract}

\section{Introduction}
\label{sec:Intro}
Eigenvalue problems in the context of vibroacoustics often depend on several parameters.
In this work, we consider a geometry and material dependent violin bridge.
For a fast and reliable evaluation in the real-time and multi-query context, reduced basis methods have proven to be a powerful tool.

For a comprehensive review on reduced basis methods, see, e.\,g.,~\cite{QaMaNe15,RoHuPa08}
or~\cite[Chapter~19]{Qua14} and the references therein.
The methodology has been applied successfully to many different problem classes, among others Stokes problems~\cite{IQRV14,LoMaRo06,RoHuMa13,RoVe07}, variational inequalities~\cite{GU14,HaSaWo12} and linear elasticity~\cite{milani:08}.
Recently, reduced basis methods for parameterized elliptic eigenvalue problems (\muEVP s) gained attention.
 
Early work on a residual a posteriori estimator  for the first eigenvalue can be found in~\cite{MaMaOlPaRo00} and has been generalized in ~\cite{Pau07a,Pau07,Pau08} to the case of several single eigenvalues with special focus to applications in electronic structure problems in solids.
Furthermore, the very simple and special case of a single eigenvalue where only the  mass matrix and not the stiffness matrix of a generalized eigenvalue problem is parameter dependent has been discussed in~\cite{FuMaPaVe15}.
Alternatively to the classical reduced basis approach, component based reduction strategies are considered in~\cite{VaHuKnNgPa14}. Here, we follow the ideas of~\cite{horger:16} where rigorous bounds  in the case of multi-query  and multiple eigenvalues are given.

More precisely, a single reduced basis is built for all eigenvalues of interest, based on a greedy selection using an online-offline decomposable error estimator.

The eigenvalues of  a violin bridge play a crucial role in transmitting the vibration of the strings to the violin body and hence influence the sound of the instrument, see~\cite{FlRo98,woodhouse:05}.  Due to the complicated curved domain and improved eigenvalue approximations compared to finite element methods, see~\cite{HuEvRe14}, we consider an isogeometric discretization. Flexibility for the tensor product spline spaces are gained by a weak domain decomposition of the non-convex domain. 

Isogeometric analysis, introduced in 2005 by  Hughes et al. in~\cite{hughes:05}, is a family of  methods that uses B-splines and non-uniform rational B-splines (NURBS) as basis functions to construct numerical approximations of partial differential equations, see also~\cite{beirao:14,hughes:09}. Mortar  methods are a popular tool for the weak coupling of non-matching meshes, originally introduced for spectral and finite element methods~\cite{ben_belgacem:99,ben_belgacem:97,bernardi:94}.
An early contribution to isogeometric elements in combination with domain decomposition techniques can be found in~\cite{hesch:12}. A rigorous mathematical analysis of the a priori error in combination with uniform stability results for different Lagrange multiplier spaces is given in~\cite{brivadis:15} and applications of isogeometric mortar methods can be found in~\cite{dittmann:14, dornisch:14, seitz:16}.

Mortar formulations are quite often formulated as an indefinite saddle point-problem. The additional degrees of freedom for the Lagrange multiplier  as well as the need for a uniform inf-sup condition to achieve stability make mortar methods, in general, more challenging than simple conforming approaches. Theoretically, the Lagrange multiplier can be eliminated. However this often results in a global process and is not carried out directly. 
While this concerns the detailed solution, the reduced basis can be purely based on a primal space and results in a non-conforming but positive definite approach.
Then the saddle point structure becomes redundant, and we gain the efficiency of a positive definite reduced system.

This article is structured as follows. In Section~\ref{sec:ProblemSetting}, we describe the geometric setup and the isogeometric mortar discretization for the violin bridge.
The model order reduction is introduced in Section~\ref{sec:ReducedBasis}, including a modified error estimator.
Numerical results illustrating the accuracy and flexibility of the presented approach are given in Section~\ref{sec:Numeric}. 
\section{Problem setting}
\label{sec:ProblemSetting}
\begin{figure}
\begin{minipage}[b]{.4\textwidth}
\includegraphics[ height =  8em ]{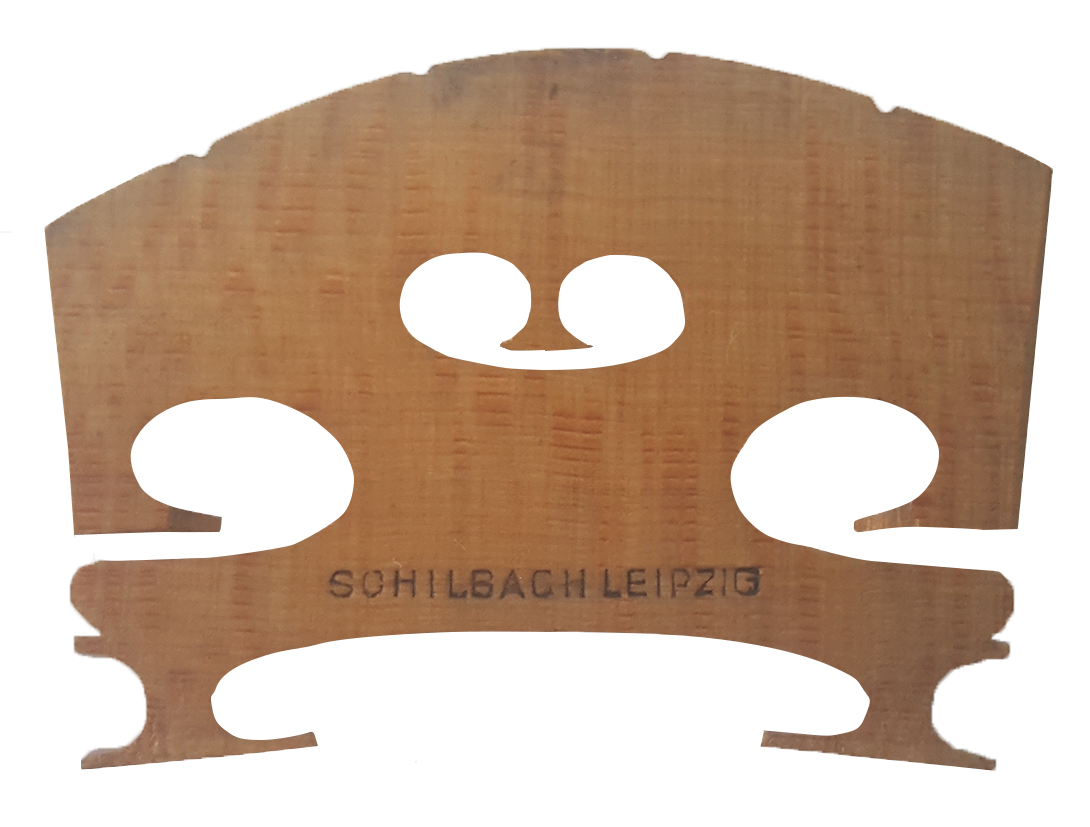}

\begin{tikzpicture}[overlay]

\tikzstyle{style1} = [black]
\tikzstyle{style2} = [red, very thick, dashed]
\tikzstyle{style3} = [blue, very thick, dash pattern=on \pgflinewidth off 1pt]

\def\x{.1}
\def\y{.4}
\def\a{1.0};
\def\b{1.0};

\coordinate (0) at (\x, \y);
\coordinate (1) at (\x+\b, \y);
\coordinate (2) at (\x, \y+\a);

\draw[thick,->]  (0) -- (1);
\draw[thick,->]  (0) -- (2);
\def\lnb{0.4};
\def\lna{1};
\def\lnc{0.15};
\node[style1] at ($(\x+.5*\a,\y-\lnc)$) {${\bf x}$};
\node[style1] at ($(\x-\lnc,\y+.5*\b)$) {${\bf y}$};
\end{tikzpicture}
\caption{Example of a violin bridge.}
\label{fig:violin_picture}
\end{minipage}\hfill
\begin{minipage}[b]{.553\textwidth}
\includegraphics[ height =  9.5em ]{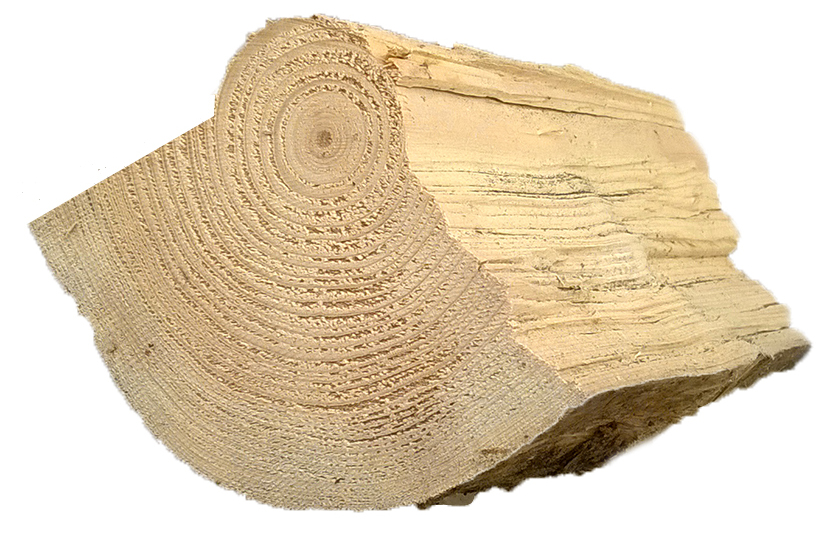}

\begin{tikzpicture}[overlay]

\tikzstyle{style1} = [black]
\tikzstyle{style2} = [red, very thick, dashed]
\tikzstyle{style3} = [blue, very thick, dash pattern=on \pgflinewidth off 1pt]

\def\x{3.75};
\def\y{2};

\def\a{3.2};
\def\b{1.6};

\def\c{.25};
\def\d{2.6};

\def\e{2};
\def\f{-2.5};

\def\lenx{.3};
\def\leny{.25};
\def\lenz{.25};

\coordinate (0) at (\x, \y);
\coordinate (1) at (\x+\lenx*\d, \y+\lenx*\c);
\coordinate (2) at (\x+\leny*\b, \y+\leny*\a);
\coordinate (3) at (\x+\lenz*\e, \y+\lenz*\f);

\draw[thick,->]  (0) -- (1);
\draw[thick,->]  (0) -- (2);
\draw[thick,->]  (0) -- (3);
\def\lnb{0.4};
\def\lna{1};
\def\lnc{1.1};
\node[style1] at ($(\x+.5*\lenx*\d,\y+.5*\lenx*\c+.1)$) {${\bf y}$};
\node[style1] at ($(\x+.5*\leny*\b-.1,\y+.5*\leny*\a)$) {${\bf z}$};
\node[style1] at ($(\x+.5*\lenz*\e-.15,\y+.5*\lenz*\f-.05)$) {${\bf x}$};
\end{tikzpicture}
\caption{Illustration of the orthotropic structure of wood.}
\label{fig:orth_wood}
\end{minipage}
\end{figure}

In vibroacoustical applications, often complicated curved domains are of special interest. Besides large constructions, such as bridges and buildings, also music instruments are investigated. An important part of the violin is a wooden violin bridge, depicted in Figure~\ref{fig:violin_picture}. 
For this geometry, we consider the eigenvalue problem of elasticity
\[
-\divergence \sigma(u) = \ew \rho u,
\]
where $\rho>0$ is the mass density, and $\sigma(u)$ depends on the material law of the structure under consideration. In our case, linear orthotropic materials are appropriate since as depicted in Figure~\ref{fig:orth_wood} wood consists of three different axes and only small deformations are considered. Note that besides the cylindrical structure of a  tree trunk, we consider Cartesian coordinates due to the comparably small size of the bridge.

\subsection{Orthotropic material law} \label{subsec:orthotrop}
The three axes are given by the fiber direction $y$, the in plane orthogonal direction $z$ and the radial direction $x$. By Hooke's law, the stress strain relation can be stated in its usual form as $\sigma(u) = \mathbb{C} \varepsilon(u)$ with $\varepsilon(u) = (\nabla u + \nabla u ^\top)/2$. Due to the alignment of the coordinate system with the orthotropic structure, the stiffness tensor is given as
\begin{equation} \label{eq:stiffness_tensor}
 \mathbb{C}=\begin{pmatrix}
 A_{11} & A_{12} & A_{13}  & 0&0 & 0 \\
 A_{21} & A_{22} & A_{23}  &0 &0 & 0 \\
 A_{31} & A_{32} & A_{33}  &0 &0 & 0 \\
 0 & 0 &  0 & G_{yz} &0 & 0 \\
 0 & 0 & 0  &0 &G_{zx} &0  \\
 0 & 0 &  0 &0 & 0& G_{xy}
 \end{pmatrix},
\end{equation}
with the shear moduli $G_{xy}, G_{yz}, G_{zx}$ and the entries $A_{ij}$ depending on the elastic moduli $E_x, E_y, E_z$ and the Poisson's ratios $\nu_{xy}, \nu_{yz}, \nu_{zx}$. The exact formula for $A_{ij}$ can be found in~\cite[Chapter 2.4]{rand:07}.

Some important differences compared to isotropic material laws are worth pointing out. 
While in the isotropic case, all Poisson's ratios share the same value, for orthotropic materials they represent three independent material parameters.
The only relation between the ratios is $\nu_{ij} E_{j} = \nu_{ji} E_{i}$. Also the possible range  of the material parameters, i.e., $-1<\nu<1/2$  for the isotropic case, is different. 
$1>\nu_{yz}^2{E_z}/{E_y}+\nu_{xy}^2 {E_y}/{E_x}+ 2\nu_{xy}\nu_{yz}\nu_{zx} {E_z}/{E_x}+\nu_{zx}^2{E_z}/{E_x}$ and $E_{x}/E_{y}>\nu_{xy}^2$ are necessary for a positive definite stiffness tensor and hence a coercive energy functional. Note that Poisson's ratios larger than $1/2$ are permitted, but this does not imply unphysical behavior as in the isotropic case, see, e.g.,~\cite{ranz:07}. 
The conditions $E_i, G_{ij} >0$ hold both in the isotropic and orthotropic case. 

\begin{figure}
\centering
\includegraphics[ width = .7 \textwidth ]{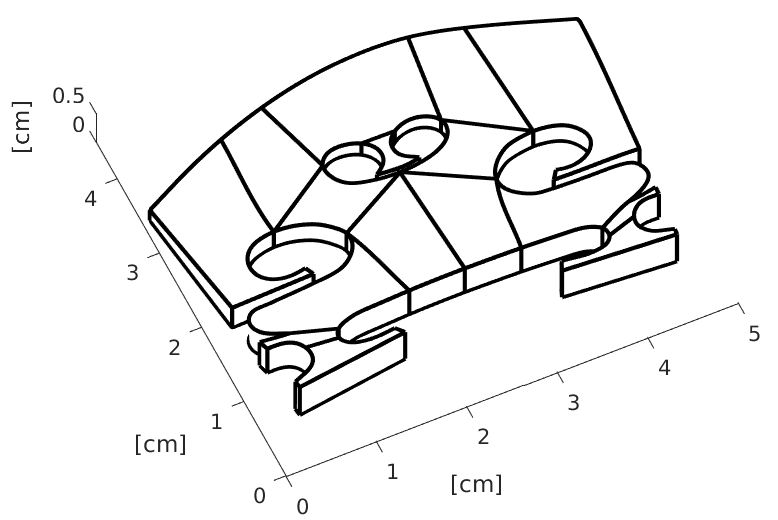}
\caption{Decomposition of the three-dimensional geometry into 16 patches and 16 interfaces.}
\label{fig:bridge_decomp_3d}
\end{figure}

The curved domain of the violin bridge can be very precisely described by a spline volume. Since it is not suitable for a single-patch description, we decompose it into 16 three-dimensional spline patches shown in Figure~\ref{fig:bridge_decomp_3d}.  While the description of the geometry could also be done with fewer patches, the number of 16 patches $\Omega_i$ gives us regular geometry mappings and a reasonable flexibility of the individual meshes. 
The decomposed geometry is solved using an equal-order isogeometric mortar method as described in~\cite{brivadis:15}. A trivariate B-spline space $V_i$ is considered on each patch $\Omega_i$. The broken ansatz space $V_h = \prod_i V_i$ is weakly coupled on each of the 16 interfaces. For each interface $\gamma_k$ the two adjacent domains are labeled as one slave and one master domain (i.e. $\gamma_k = \partial \Omega_s \cap \partial \Omega_m$) and the coupling space $M_k$ is set as the trace space  of the spline spaces on the slave domain and $M_h = \prod_k M_k$. Several crosspoints are present in the decomposition, where  an appropriate local degree reduction is performed, as described in~\cite[Section 4.3]{brivadis:15}, to guarantee stability.

We use the standard bilinear forms for mortar techniques in linear elasticity
\[
a(u,v) = \sum_i \int_{\Omega_i} \sigma(u):\varepsilon(v), \quad 
m(u,v) = \sum_i \int_{\Omega_i} \rho u v, \quad
b(v, \lmtest) = \sum_k \int_{\gamma_k} [v]_k \lmtest,
\]
where $[v]_k = \left. v_s \right|_{\gamma_k} - \left. v_m \right|_{\gamma_k}$ denotes the jump across the interface $\gamma_k$. We note that no additional variational crime by different non-matching geometrical resolutions of $\gamma_k$ enters. The detailed eigenvalue problem is defined as $(u,\lm)\in V_h\times M_h$, $\lambda\in\R$, such that 
\begin{subequations}\label{eq:detailed_spp}
\begin{align}
a(u,v) + b(v,\lm) &= \ew m(u,v), \quad v\in V_h, \\
b(u, \lmtest) & = 0, \quad \lmtest \in M_h.
\end{align}
\end{subequations}
The Lagrange multiplier $\lm$ is an approximation for the surface tension $\sigma(u) n$ along the interfaces.

Additional to the nine material parameters $E_i$, $G_{ij}$, $\nu_{ij}$, we consider a geometry parameter $\mu_{10}$, describing the thickness of the violin bridge. Transforming the geometry to a reference domain, we can interpret the thickness parameter as one more material parameter. 

\subsection{Transforming geometrical parameters to material parameters}
 Let the parameter dependent geometry $\Omega(\mu)$ be an 
unidirectional scaling  of a reference domain $\hatOmega$, i.e., a transformation by $F(\cdot; \mu):\hatOmega\rightarrow\Omega(\mu)$, ${\bf x} = F({\hatx}; \mu) = (\widehat x,\widehat y, \mu_{10} \widehat z)$, with $\hatx = (\widehat x, \widehat y, \widehat z)\in \widehat \Omega$. Transforming the unknown displacement and rescaling it as 
$\hatu(\hatx ) = \DF u(F(\hatx;\mu))$
allows us to define a symmetric strain variable on the reference domain
\[\hatepsilon(\hatu(\hatx)) = \DF  \epsilon(u(F(\hatx;\mu))) \DF.\]
The orthotropic stiffness tensor~\eqref{eq:stiffness_tensor} is then transformed to
\begin{align*}
\widehat{\mathbb{C}}(\mu) = \begin{pmatrix}
\phantom{\mu_{10}^{-2}}A_{11} &\quad\phantom{\mu_{10}^{-2}}A_{12}&\quad \mu_{10}^{-2}A_{13} & & &\\ \phantom{\mu_{10}^{-2}}A_{21} &\quad\phantom{\mu_{10}^{-2}}A_{22}&\quad \mu_{10}^{-2}A_{23} &&&\\\mu_{10}^{-2}A_{31}&\quad  \mu_{10}^{-2}A_{32}&\quad\mu_{10}^{-4}A_{33} &&&\\ &&&&\mu_{10}^{-2} G_{yz}&&\\&&&&&\mu_{10}^{-2}G_{zx}& \\ &&&&&&G_{xy}
\end{pmatrix}.
\end{align*} 
In terms of this coordinate transformation, the eigenvalue problem in the continuous  $H^1$-setting reads, since $\det \DF = \mu_{10}^{-1}$ is constant,  as
\[\int_\hatOmega {\hatepsilon(\hatu)}~\widehat{\mathbb{C}}(\mu) ~{\hatepsilon(\hatv)} ~ \mathrm{d}\hatx =
\lambda \int_\Omega \rho~ \hatu^\top \begin{pmatrix} 1&&\\&1&\\&&\mu_{10}^{-2}\end{pmatrix} \hatv ~ \mathrm{d}\hatx.\]

In our mortar case, the  coupling conditions across the interfaces have to be transformed as well. However due to the affine mapping with respect to x, y and z,
the parameter dependency on $b(\cdot,\cdot)$ can be removed.  

Another material parameter of interest for applications is the constant mass density $\rho$. However any change in the constant $\rho$ does not influence the eigenvectors. Only the eigenvalue is rescaled, yielding a trivial parameter dependence. For this reason, the density is kept constant in the reduced basis computations and can be varied in a postprocess by rescaling the eigenvalues.

The described material parameters allow for an affine parameter dependence of the mass and the stiffness, with $Q_a = 10$, $Q_m = 2$,
 \[
 a(\cdot, \cdot;\mu) = \sum_{q=1}^{Q_a} \theta^q_a(\mu) a^q(\cdot, \cdot), \quad m(\cdot, \cdot;\mu) = \sum_{q=1}^{Q_m} \theta_m^q(\mu) m^q(\cdot, \cdot).
 \]

\section{Reduced Basis}
\label{sec:ReducedBasis}
Reduced basis methods for the simultaneous approximation of eigenvalues and eigenvectors have been analyzed in~\cite{horger:16}. Here, we apply these methods to the previously described setting. An important difference to the previous work is that we wish to approximate a saddle point problem instead of a positive definite matrix.  

Previous works on saddle point problems construct a reduced basis for both the primal and the dual space. This is necessary for example for variational inequalities or when the coupling is parameter dependent, see \cite{gerner:12, glas:15, HaSaWo12, negri:15}. To ensure the inf-sup stability of the discrete saddle point problem, supremizers can be added to the primal space, additionally increasing the size of the reduced system.

\subsection{Reduction of the saddle point problem}

Due to the parameter-independence of $b(\cdot,\cdot)$, we can 
 reformulate the detailed saddle point problem~\eqref{eq:detailed_spp} in a purely primal form posed on the constrained space $X_h = \{v \in V_h, b(v, \lmtest) = 0, \lmtest \in M_h\}$. 
Note that this formulation is not suitable for solving the detailed solution, since, in general, it is costly to construct a basis of $X_h$ and severely disturbs the sparsity of the detailed matrices. Only in the case of so-called dual Lagrange multiplier spaces, a local static condensation can be carried out and the constrained basis function do have local support. However, in the reduced basis context the constructed basis functions do automatically satisfy the weak coupling properties and thus the saddle-point problem is automatically reduced to a positive definite one.

Our reduced space is defined by 
$X_N = \{\zeta_n \in X_h, n=1,\ldots,N\}$, where the reduced basis functions $\zeta_n$ are selected as presented in~\cite{horger:16}.
Then the reduced eigenvalue problem for the first $K$ eigenpairs is given by: 
Find the eigenvalues $\lambdaRB{i}(\mu) \in \R$ and the eigenfunctions  $\uRB{i}(\mu) \in X_N$, $i=1,\ldots,K$, such that 
\begin{align*} 
a(\uRB{i}(\mu),v;\mu) &= \lambdaRB{i}(\mu) ~ m(\uRB{i}(\mu),v;\mu), \quad v\in X_N
\end{align*}

In a first step, an initial basis is built by a small POD. This basis is then enlarged by a greedy algorithm with an  asymptotically reliable error estimator. 

\subsection{Decomposition of the error estimator with a parameter-dependent mass} 
The error estimator presented in~\cite[Corollary 3.3]{horger:16} can directly be applied, but the online-offline decomposition needs to be modified. In the original setting, a parameter-independent mass was considered, so we need to additionally include the affine decomposition of the mass matrix.

The main contribution of the estimator is the residual  \[r_i(\cdot;\mu) = a(\uRB{i}(\mu), \cdot; \mu) - \lambdaRB{i}(\mu)\, m(\uRB{i}(\mu), \cdot; \mu)\] measured in the dual norm $  \aNormRefDual{r} = \sup_{v\in X_h} {r(v)}/{\aRef(v,v)}^{1/2}$, where $\aRef(u,v) := a(u,v;\widehat \mu)$, and $\widehat \mu \in \mathcal{P}$ is a reference parameter. We define $ \eRef_i(\mu)\in X_h$ by \[ a( \eRef_i(\mu), v; \mu) = r_i(v;\mu),\quad v\in X_h.\]

To adapt the online-offline decomposition, we follow~\cite{horger:16,MaMaOlPaRo00} and add additional terms corresponding to the mass components $m_q$.
Let $(\zeta_n)_{1\leq n\leq N}$ be the orthonormal basis (w.\,r.\,t.\ $m(\cdot,\cdot;\widehat \mu)$) of $X_N$ and let us define $\xi_n^q\in X_N$ and $\xi_n^{m,q}\in X_N$ by
\begin{align*}
  \aRef(\xi_n^q,v) &= a^q(\zeta_n,v),\quad v\in X_h,\;1\leq n\leq N,\;1\leq q\leq Q_a, \\
  \aRef(\xi_n^{m,q},v) &= m^q(\zeta_n,v),\quad v\in X_h,\;1\leq n\leq N,\;1\leq q\leq Q_m.
\end{align*}
In the following,
we identify the function $\uRB{i}(\mu)\in V_N$
and its vector representation w.\,r.\,t.\ the basis $(\zeta_n)_{1\leq n\leq N}$
such that
$\left(\uRB{i}(\mu)\right)_n$ denotes the $n$-th coefficient.
Then,
given a reduced eigenpair $\left(\uRB{i}(\mu),\lambdaRB{i}(\mu)\right)$, we have the error representation
\begin{equation*}
  \eRef_i(\mu) = \sum_{n=1}^N \sum_{q=1}^{Q_a} \theta_a^q(\mu) \left(\uRB{i}(\mu)\right)_n \xi_n^q
  - \lambdaRB{i}(\mu) \sum_{n=1}^N  \sum_{q=1}^{Q_m}  \theta_m^q(\mu) \left(\uRB{i}(\mu)\right)_n  \xi_n^{m,q}.
\end{equation*}
Consequently,
the main contribution of $\eta_i(\mu)$ decomposes using $ \aNormRefDual{r_i(\cdot;\mu)}^2 =  \aRef(  \eRef_i(\mu),  \eRef_i(\mu) )$, see~\cite[Section 3.3]{horger:16} for a more detailed discussion.

\section{Numerical simulation}
\label{sec:Numeric} 
\begin{table}\hspace*{-7em}
\begin{tabular}{r||c|c|c|c|c|c|c|c|c}
 & $E_x$ $[MPa]$ &  $E_y$ $[MPa]$& $E_z$ $[MPa]$
  & $G_{yz}$ $[MPa]$ & $G_{zx}$ $[MPa]$ & $G_{xy}$ $[MPa]$ 
   & $\nu_{yz}$ & $\nu_{zx}$ &$ \nu_{xy}$ \\\hline\hline &&&&&&&&& \\[-.4em]
   $\widehat{\mu}$ & $14,000$&$2,280$& $1,160$& $465$& $1,080$& $1,640$& $0.36$& $0.0429$& $0.448$\\[1em]
$\mathcal{P}_1$ & ~ \parbox[t]{4em}{$13,000$\\ -- $15,000$}~ &  
  ~\parbox[t]{3.5em}{$1,500$ \\-- $3,000$}  ~& 
   ~\parbox[t]{3.5em}{$750$ \\-- $1,500$}   ~&  
    ~\parbox[t]{3.5em}{$100$ \\--  $1,000$} ~& 
     ~\parbox[t]{3,5em}{$500$ \\-- $1,500$} ~& 
      ~\parbox[t]{4em}{$1,000$ \\-- $2,000$}~ & 
       ~\parbox[t]{2em}{$0.3$ \\-- $0.4$}   ~& 
        ~\parbox[t]{3em}{$0.03$ \\-- $0.06$}~&
         ~\parbox[t]{2em}{ $0.4$\\ -- $0.5$}~
         \\[2em]
$\mathcal{P}_2$ & ~ \parbox[t]{4em}{$1,000$\\ -- $20,000$}~ &  
  ~\parbox[t]{3.5em}{$100$ \\-- $5,000$}  ~& 
   ~\parbox[t]{3.5em}{$100$ \\-- $2,000$}   ~&  
    ~\parbox[t]{3.5em}{$10$ \\--  $5,000$} ~& 
     ~\parbox[t]{3,5em}{$100$ \\-- $2,500$} ~& 
      ~\parbox[t]{4em}{$100$ \\-- $5,000$}~ & 
       ~\parbox[t]{2em}{$0.1$ \\-- $0.5$}   ~& 
        ~\parbox[t]{3em}{$0.01$ \\-- $0.1$}~&
         ~\parbox[t]{2em}{ $0.3$\\ -- $0.5$}~
\end{tabular}
\caption{Reference parameter and considered parameter ranges.}
\label{tab:param_values}
\end{table}
\begin{figure}
\centering
\includegraphics[ width = .4 \textwidth ]{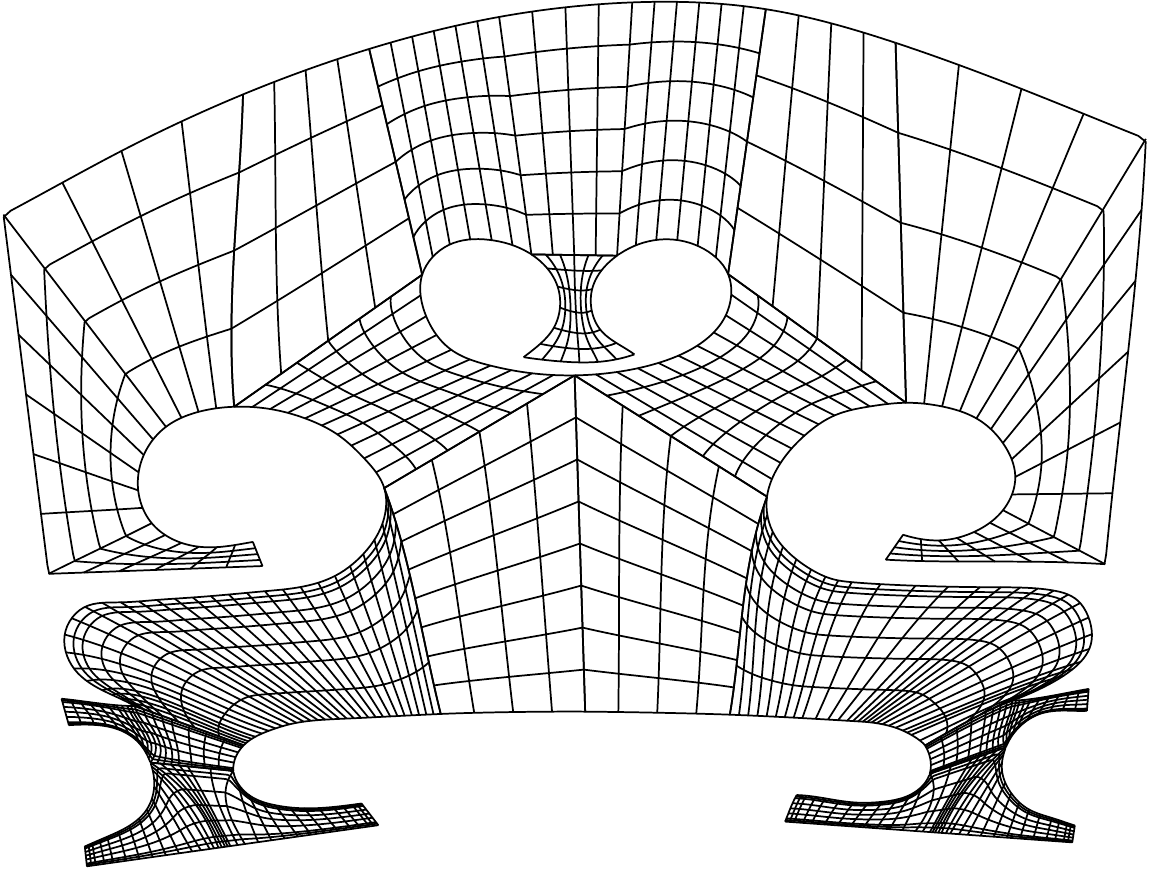}
\caption{Non-matching isogeometric mesh of the violin bridge.}
\label{fig:mesh_bridge}
\end{figure}
In this section, the performance of the proposed algorithm is illustrated by numerical examples.
The detailed computations were performed using geoPDEs~\cite{geopdes:11}, a Matlab toolbox for isogeometric analysis, the reduced computations are based on RBmatlab~\cite{DrHaKaOh12}.

For the detailed problem, we use an anisotropic discretization. In plane, we use  splines of degree $p=3$ on the non-matching mesh shown in Figure~\ref{fig:mesh_bridge}. The mesh has been adapted locally to better resolve possible corner singularities of the solution. In the $z$-direction a single element of degree $p=4$ is used. 
The resulting equation system has $45,960$ degrees of freedom for the displacement whereas the surface traction on the interfaces is approximated by $2,025$ degrees of freedom.

We consider the ten parameters, described in Section~\ref{sec:ProblemSetting}, $\mu = (\mu_1,\ldots,\mu_{10})$ with the elastic modulii $\mu_1 = E_x$, $\mu_2 = E_y$, $\mu_3 = E_z$, the shear modulii $\mu_4 = G_{yz}$, $\mu_5 = G_{xz}$, $\mu_6 = G_{xy}$, Poisson's ratios $\mu_7 = \nu_{yz}$, $\mu_8 = \nu_{xz}$, $\mu_9 = \nu_{xy}$ and the scaling of the thickness $\mu_{10}$.

The considered parameter values were chosen according to real parameter data given in~\cite[Table 7-1]{ranz:07}. We consider two 
different scenarios. In the first setting, we fix the wood type and take into account only natural variations,
see~\cite[Section 7.10]{ranz:07}. To capture the sensitivity of the violin bridge,
one can chose a rather small parameter range around a reference parameter. We chose the reference data of {\it fagus sylvatica}, the common beech, as given in Table~\ref{tab:param_values}, as well as the parameter range $\mathcal{P}_1$. 
 The mass density is fixed in all cases as $720kg/m^3$.

\begin{figure}
\centering
\includegraphics[ width = .6 \textwidth ]{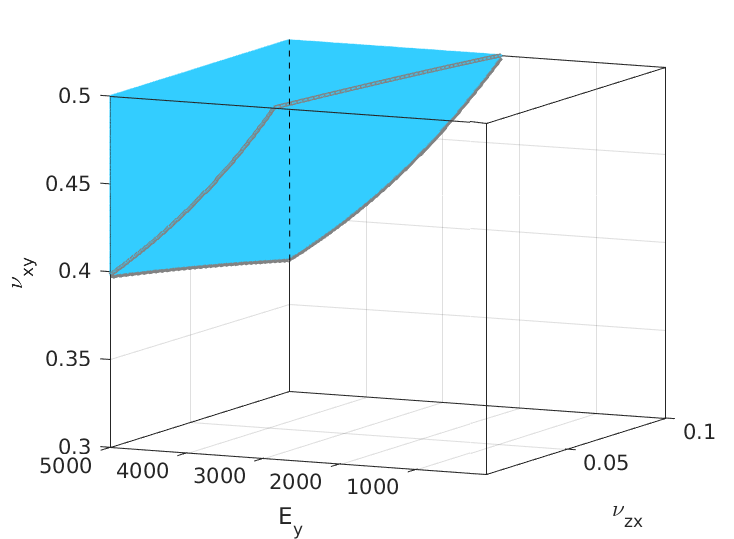}
\caption{Illustration of non-admissible parameter values in a lower-dimensional excerpt of $\mathcal{P}_2$, varying $\nu_{zx}\in (0.01, 0.1), \nu_{xy}\in (0.3,0.5), E_y\in(100,5000)$ and fixing $E_x = 1000, E_z = 2000$ and $\nu_{yz}=0.5$. }
\label{fig:parameter_admissibility}
\end{figure}

In our second test setting, we also consider different wood types.
Hence we also consider a larger parameter set, including the parameters for several types of wood. Based on a selection of some wood types, we chose the parameter range $\mathcal{P}_2$, see Table~\ref{tab:param_values}. 
We note, that not all parameters in this large range are admissible for the orthotropic elasticity as they do not fulfill the conditions for the positive definiteness of the elastic tensor, stated in Section~\ref{subsec:orthotrop}. 
Thus, we constrain the tensorial parameter space by
\[
 1-\nu_{yz}^2{E_z}/{E_y}+\nu_{xy}^2 {E_y}/{E_x}+ 2\nu_{xy}\nu_{yz}\nu_{zx} {E_z}/{E_x}+\nu_{zx}^2{E_z}/{E_x} \geq c_0, \]as well as $E_{x}/E_{y}-\nu_{xy}^2\geq c_1
$ 
where the tolerances $c_0$ and $c_1$ were chosen according to several wood types.
Exemplary, in Figure~\ref{fig:parameter_admissibility} we depict an lower-dimensional excerpt of $\mathcal{P}_2$ which includes non-admissible parameter values.

\begin{table}[htb]
\centering
\begin{tabular}{|r||c|c|c||c|c|} \hline
Eigenvalue & $\mu_{10}=0.5$ & $\mu_{10}=1.0$ & $\mu_{10}=2.0$ &~ratio $0.5$/$1.0$ &ratio $1.0/2.0$ \\ \hline \hline
$ 1$ & $ 0.4057$  &    $1.3238 $   &   $ 3.6954$  &  ~$0.3065$  & $0.3582$\\
$ 2$ & $1.1613$  &    $ 3.8870$   &   $10.8071 $  &   ~$0.2988$ & $0.3597$	\\
$ 3$ & $ 4.4096$  &    $ 12.9562 $   &   $ 26.5621 $  &  ~$0.3403$  &$0.4878$	\\
$ 4$ & $6.1371$  &    $19.3254 $   &   $30.0050 $  &   ~$0.3176$ &$0.6441$	\\
$ 5$ & $ 13.5564$  &    $ 27.3642$   &   $ 53.2657 $  &  ~$0.4954$  &$0.5137$	\\
$ 6$ & $19.2229 $  &    $46.2521 $   &   $ 93.9939$	  & ~$0.4156$   & $0.4921$ \\
$ 7$ & $27.6118 $  &    $65.0940 $   &   $ 111.6075$  &  ~$0.4242$  &$0.5832$	\\
$ 8$ & $39.3674 $  &    $96.8069 $   &   $  129.3406$   &  ~$0.4067$  & $0.7485$	\\
$ 9$ & $57.8266 $  &    $ 107.6749$   &   $ 189.6090$  &   ~$0.5370$ & $0.5679$	\\
$10$ & $ 68.0131$  &    $  130.8876$  &    $241.7695 $  &  ~$0.5196$  & $0.5414$ \\ \hline
\end{tabular}
\caption{The 10 smallest eigenvalues for different thickness parameters, with the other parameters fixed to the reference value.  }
\label{tab:param_thickness}
\end{table}

\begin{figure}[htbp!]
\begin{center}
\begin{minipage}{.32\textwidth}
\begin{center}
\includegraphics[trim={680px 100px 680px 150px},clip,width=\textwidth]{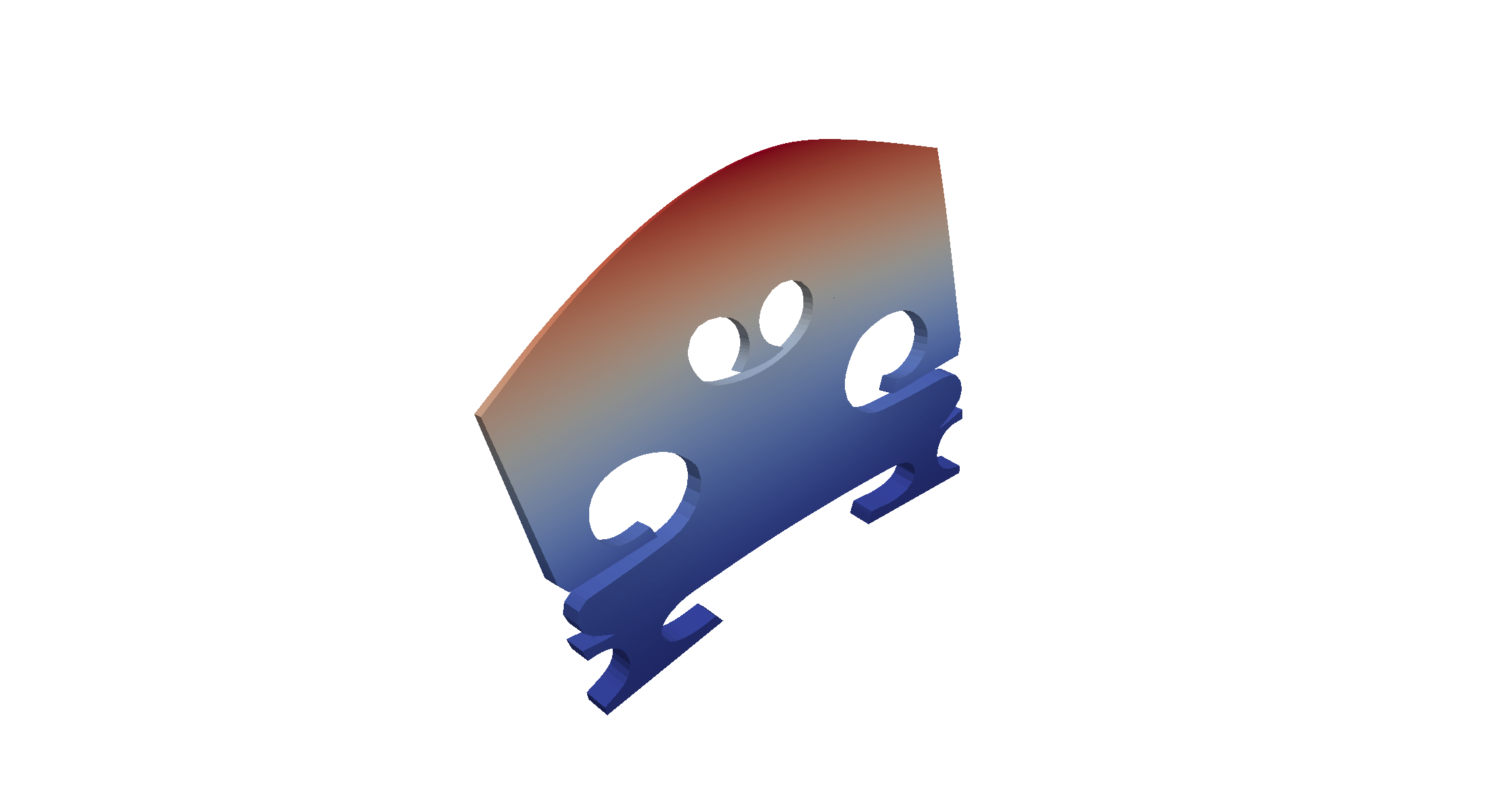}\\
first eigenvalue: $0.4057$
\end{center}
\end{minipage}
\begin{minipage}{.32\textwidth}
\begin{center}
\includegraphics[trim={680px 100px 680px 150px},clip,width=\textwidth]{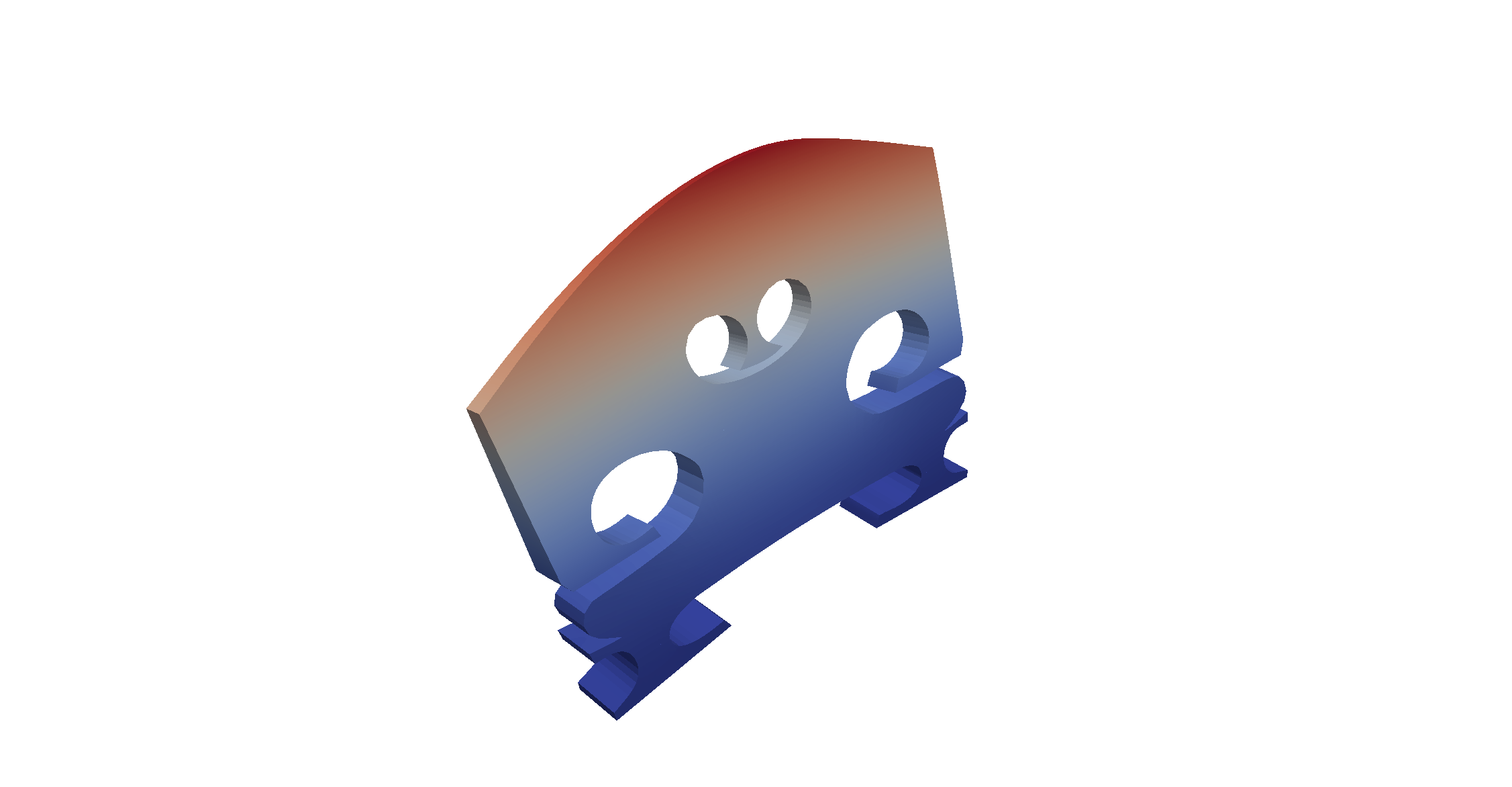}\\
first eigenvalue: $1.1613$
\end{center}
\end{minipage}
\begin{minipage}{.32\textwidth}
\begin{center}
\includegraphics[trim={680px 100px 680px 150px},clip,width=\textwidth]{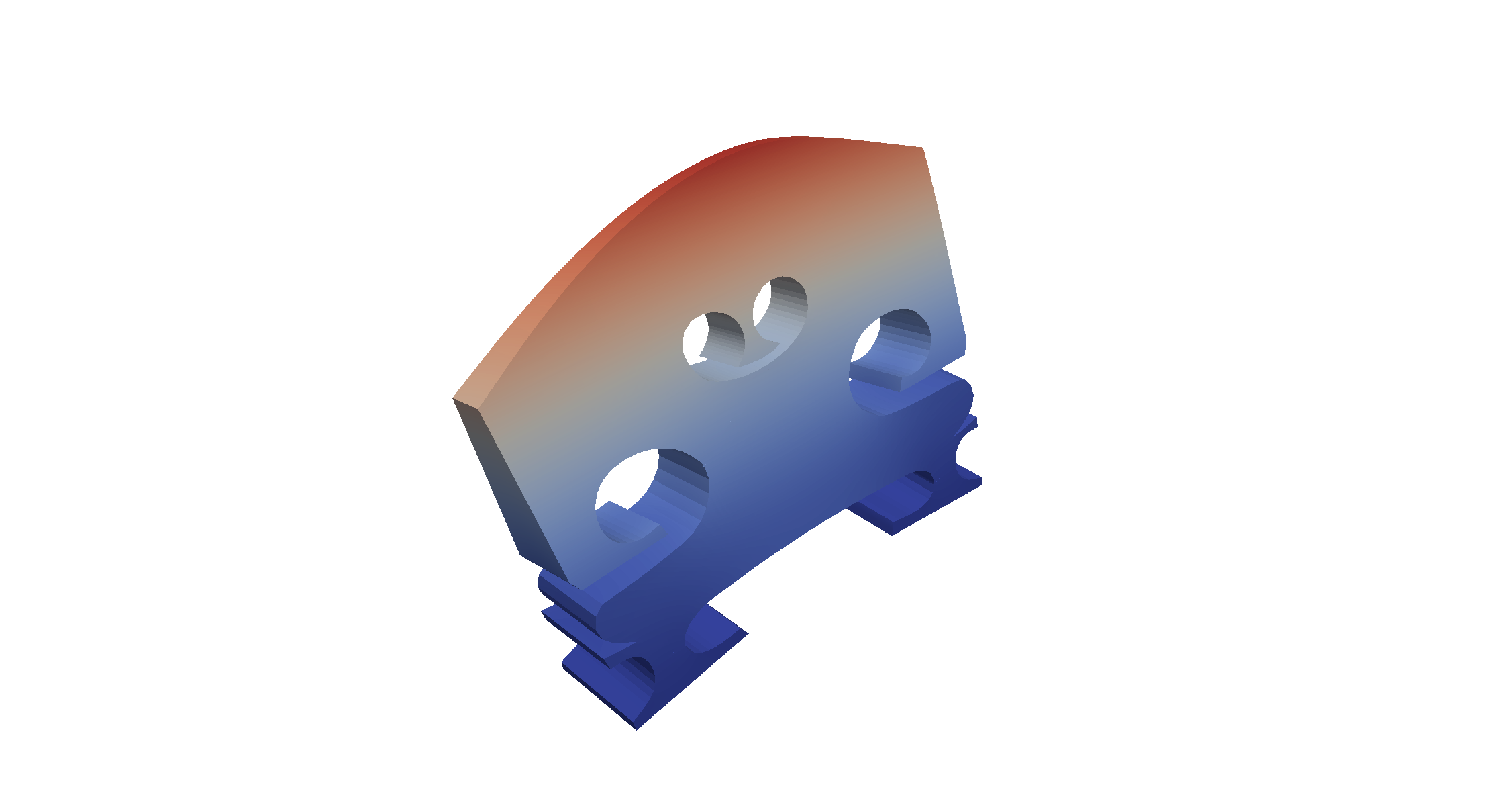}\\
first eigenvalue: $3.6954$
\end{center}
\end{minipage}\\
\begin{minipage}{.32\textwidth}
\begin{center}
\includegraphics[trim={680px 100px 680px 150px},clip,width=\textwidth]{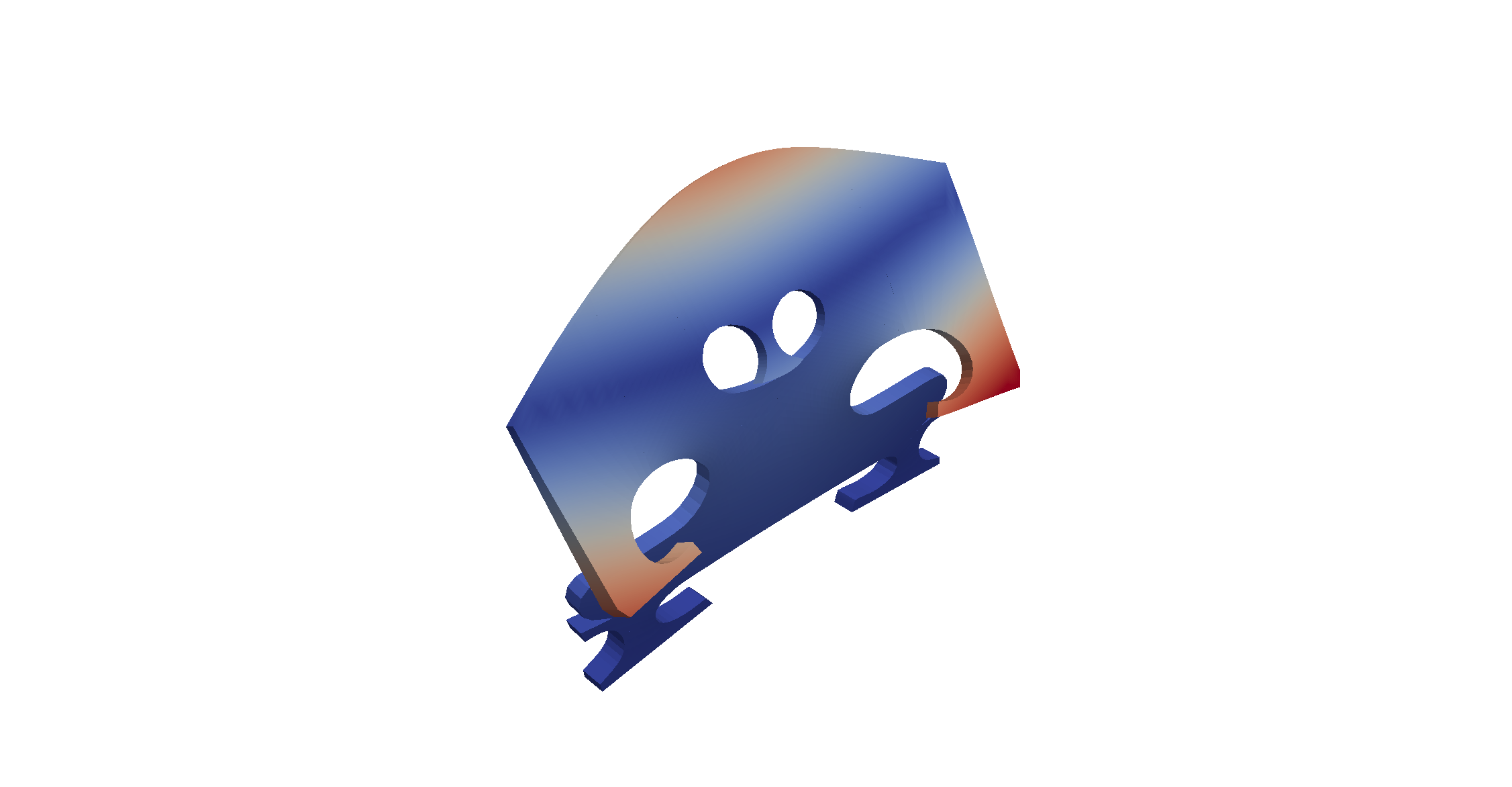}\\
third eigenvalue: $4.4096$
\end{center}
\end{minipage}
\begin{minipage}{.32\textwidth}
\begin{center}
\includegraphics[trim={680px 100px 680px 150px},clip,width=\textwidth]{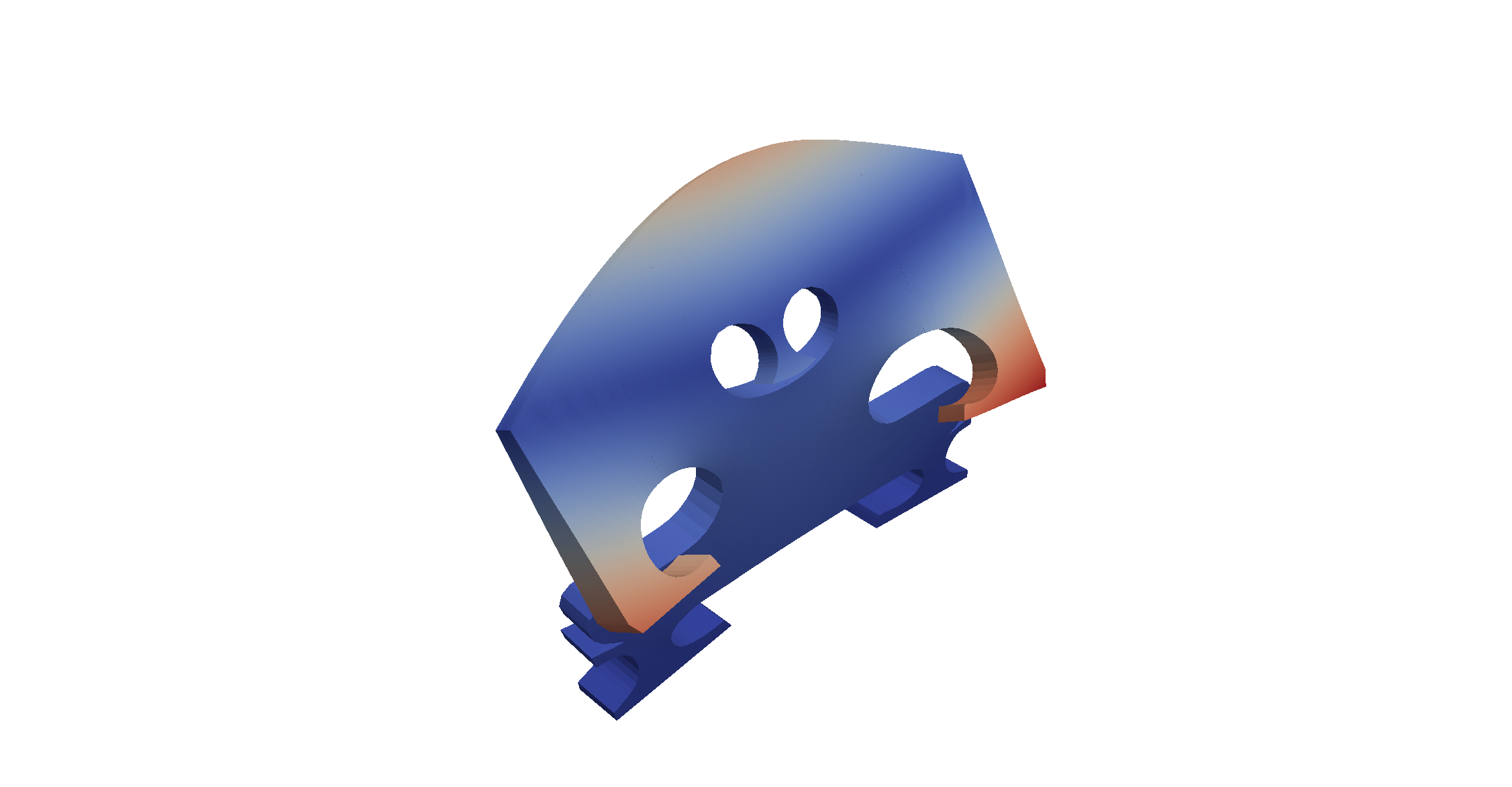}\\
third eigenvalue: $12.9562$
\end{center}
\end{minipage}
\begin{minipage}{.32\textwidth}
\begin{center}
\includegraphics[trim={680px 100px 680px 150px},clip,width=\textwidth]{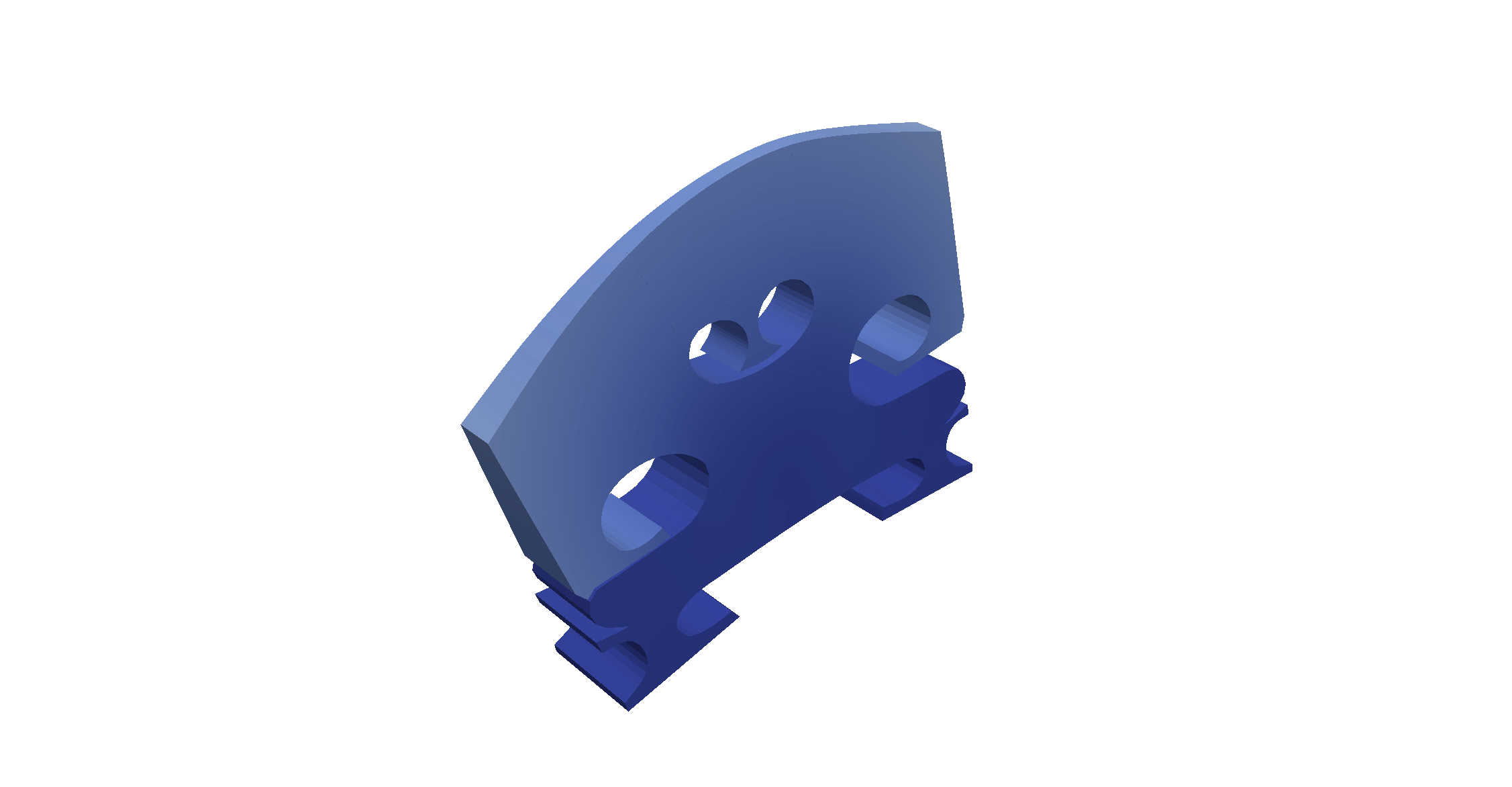}\\
third eigenvalue: $26.5621$
\end{center}
\end{minipage}\\
\begin{minipage}{.32\textwidth}
\begin{center}
\includegraphics[trim={680px 100px 680px 150px},clip,width=\textwidth]{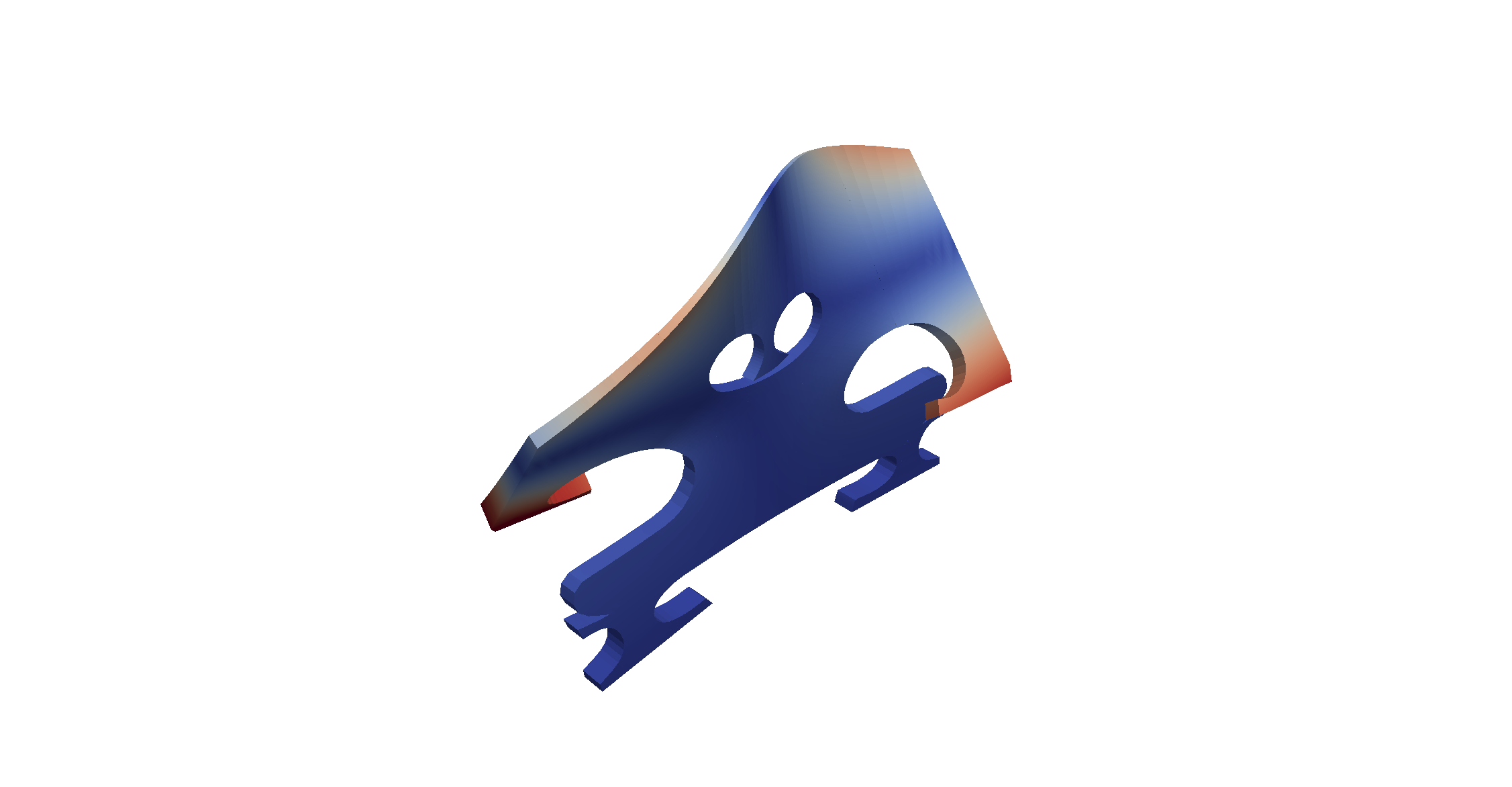}\\
fourth eigenvalue: $6.1371$
\end{center}
\end{minipage}
\begin{minipage}{.32\textwidth}
\begin{center}
\includegraphics[trim={680px 100px 680px 150px},clip,width=\textwidth]{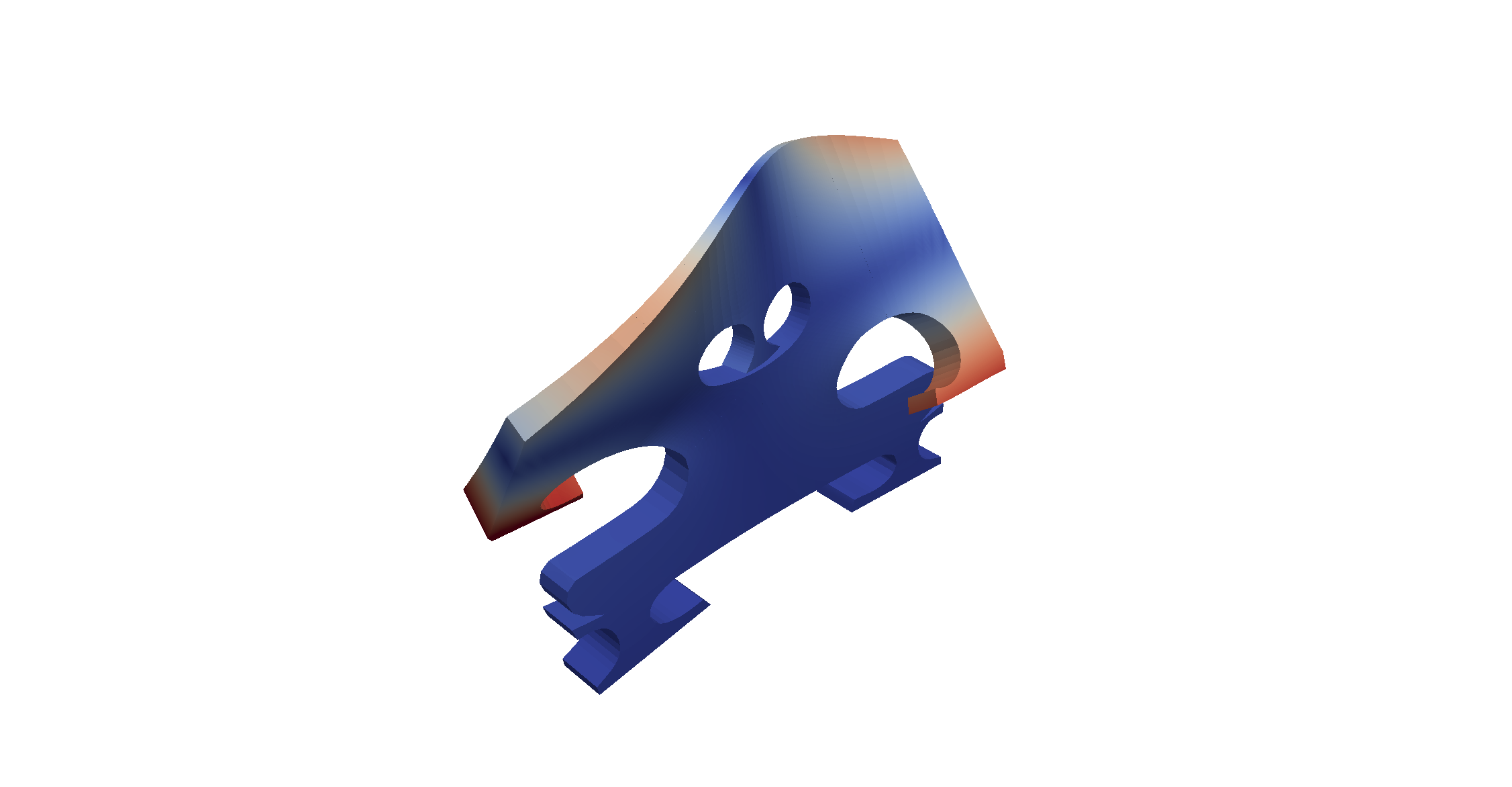}\\
fourth eigenvalue: $19.3254$
\end{center}
\end{minipage}
\begin{minipage}{.32\textwidth}
\begin{center}
\includegraphics[trim={680px 100px 680px 150px},clip,width=\textwidth]{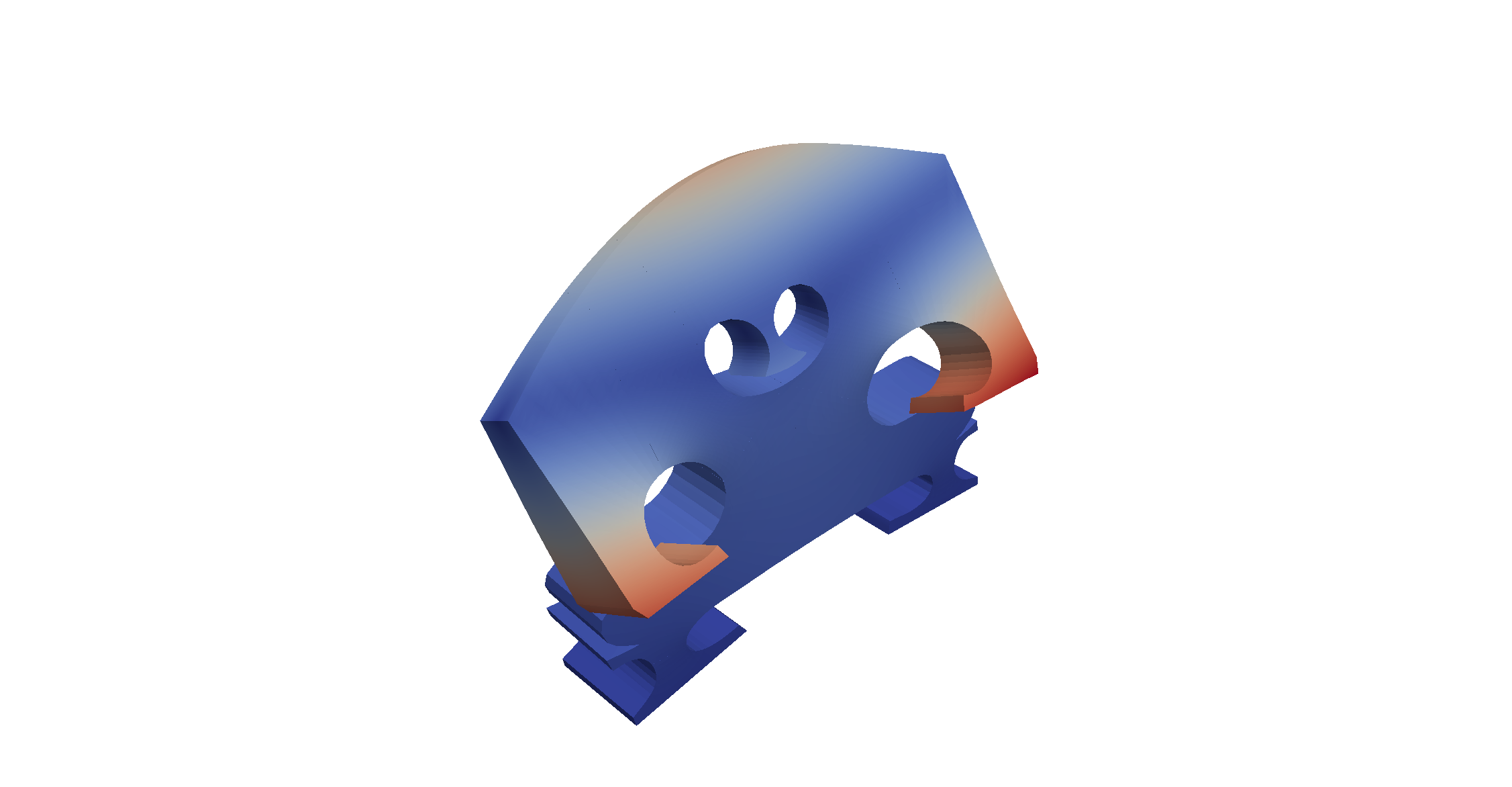}\\
fourth eigenvalue: $30.005$
\end{center}
\end{minipage}
\end{center}
\caption{Influence of the thickness of the bridge on several eigenfunctions. }
\label{fig:param_thickness}
\end{figure}

First, we consider the effect of the varying thickness parameter on the solution of our model problem. In Table~\ref{tab:param_thickness} the first eigenvalues are listed for different values of the thickness, where we observe a notable and nonlinear parameter dependency. A selection of the corresponding eigenfunctions is depicted in Figure~\ref{fig:param_thickness}, where the strong influence becomes even more evident, since in some cases the shape of the eigenmode changes when varying the thickness. 

In the following reduced basis tests, the relative error values are computed as the mean value over a large amount of random parameters. The $L^2$-error of the normed eigenfunctions is evaluated as the residual of the $L^2$-projection onto the corresponding detailed eigenspace. This takes into account possible multiple eigenvalues and the invariance with respect to a scaling by $(-1)$.

\begin{figure}[htb]
\centering
\includegraphics[width=.32\textwidth]{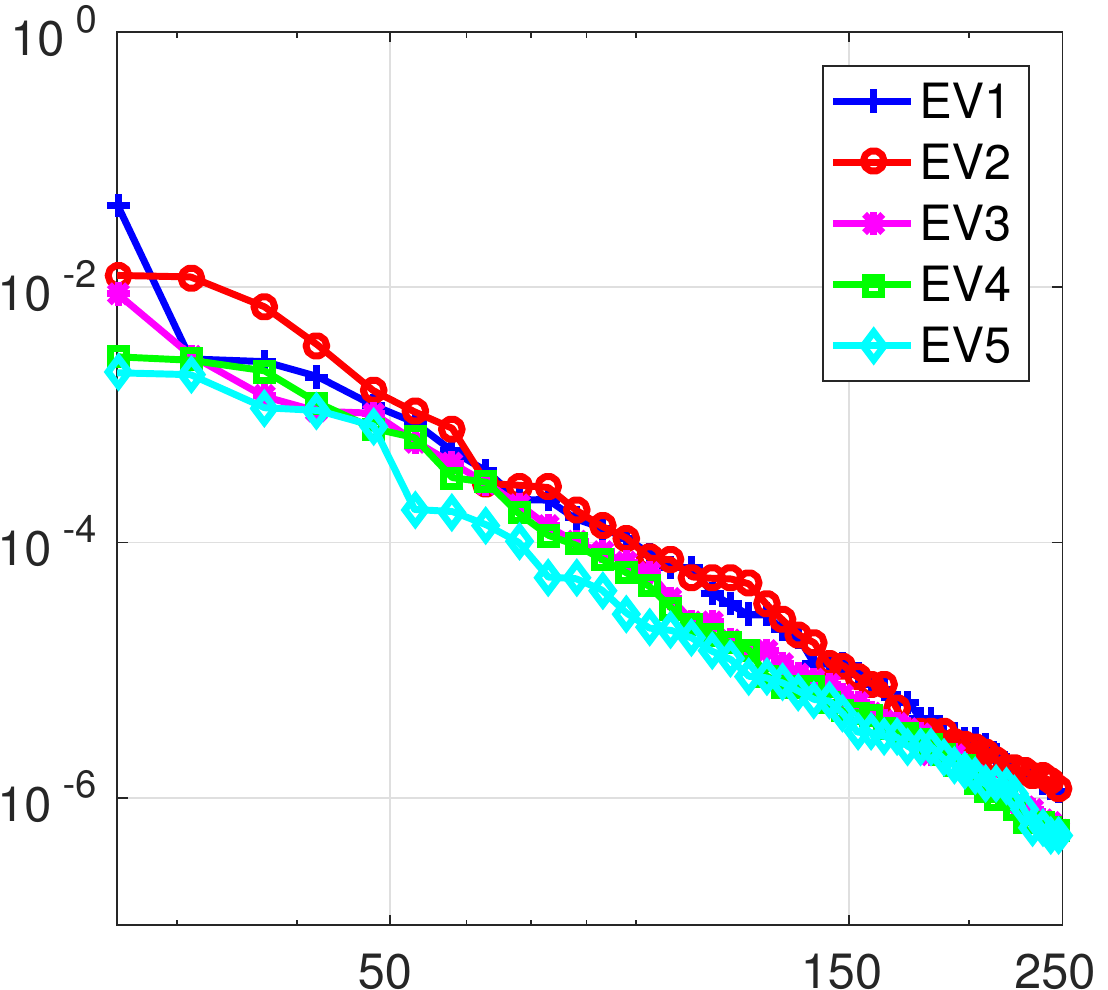}
\includegraphics[width=.32\textwidth]{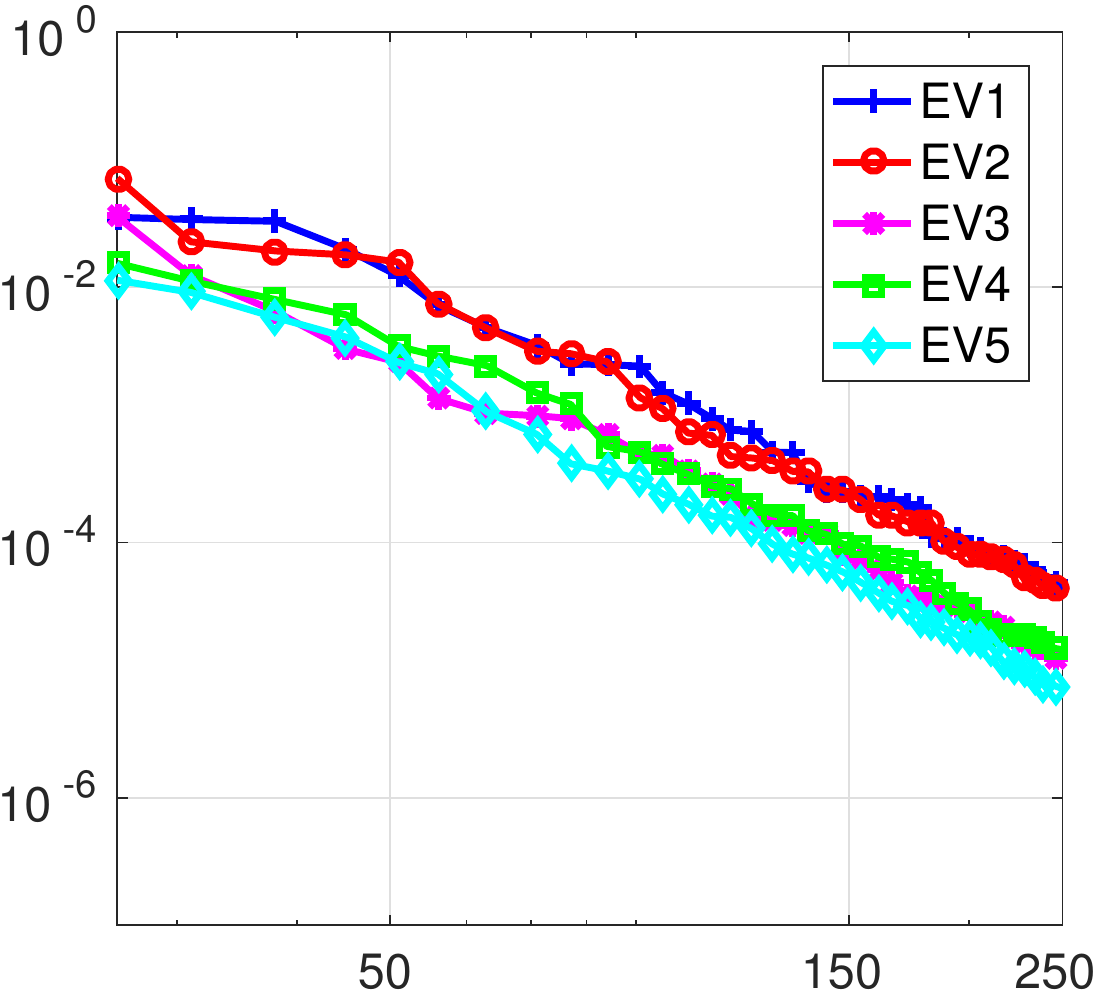}
\includegraphics[width=.32\textwidth]{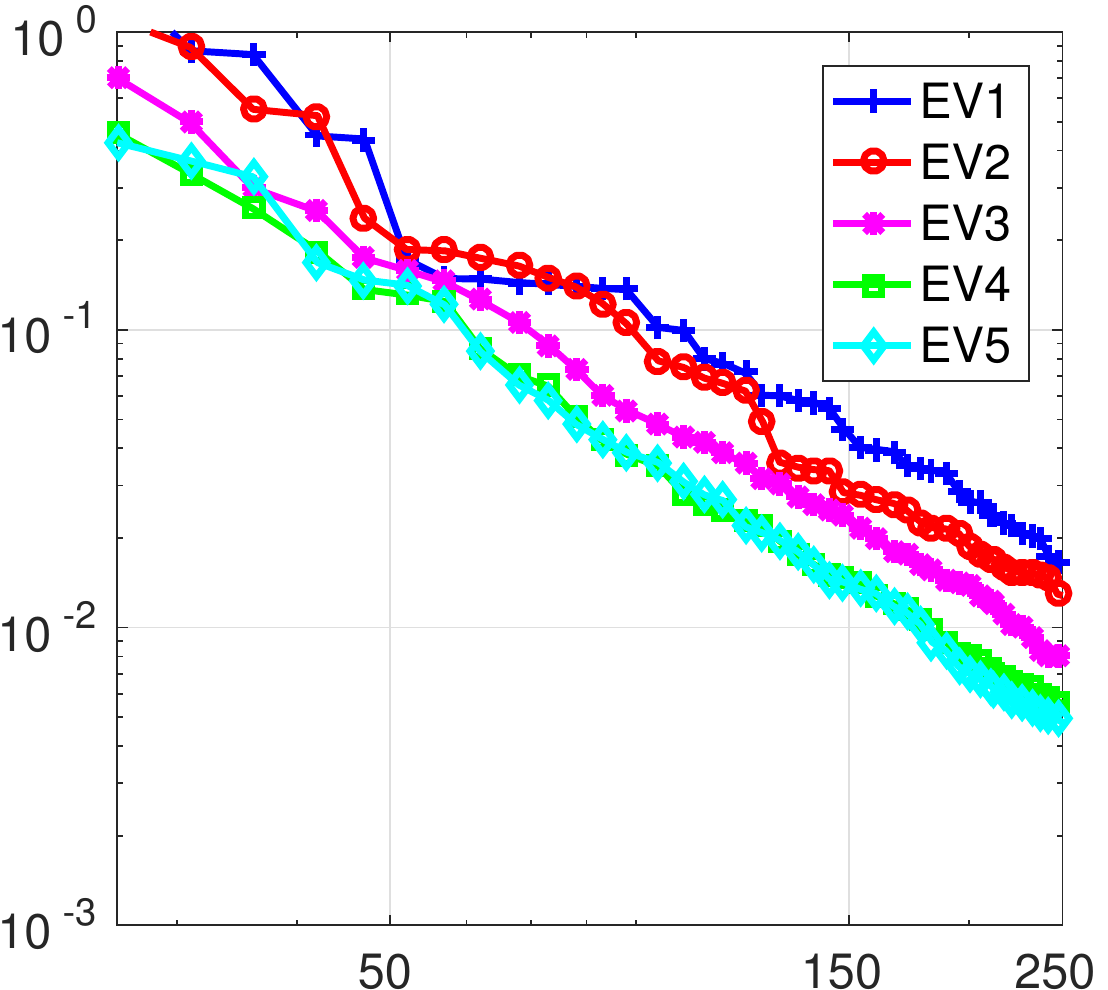}\\
\includegraphics[width=.32\textwidth]{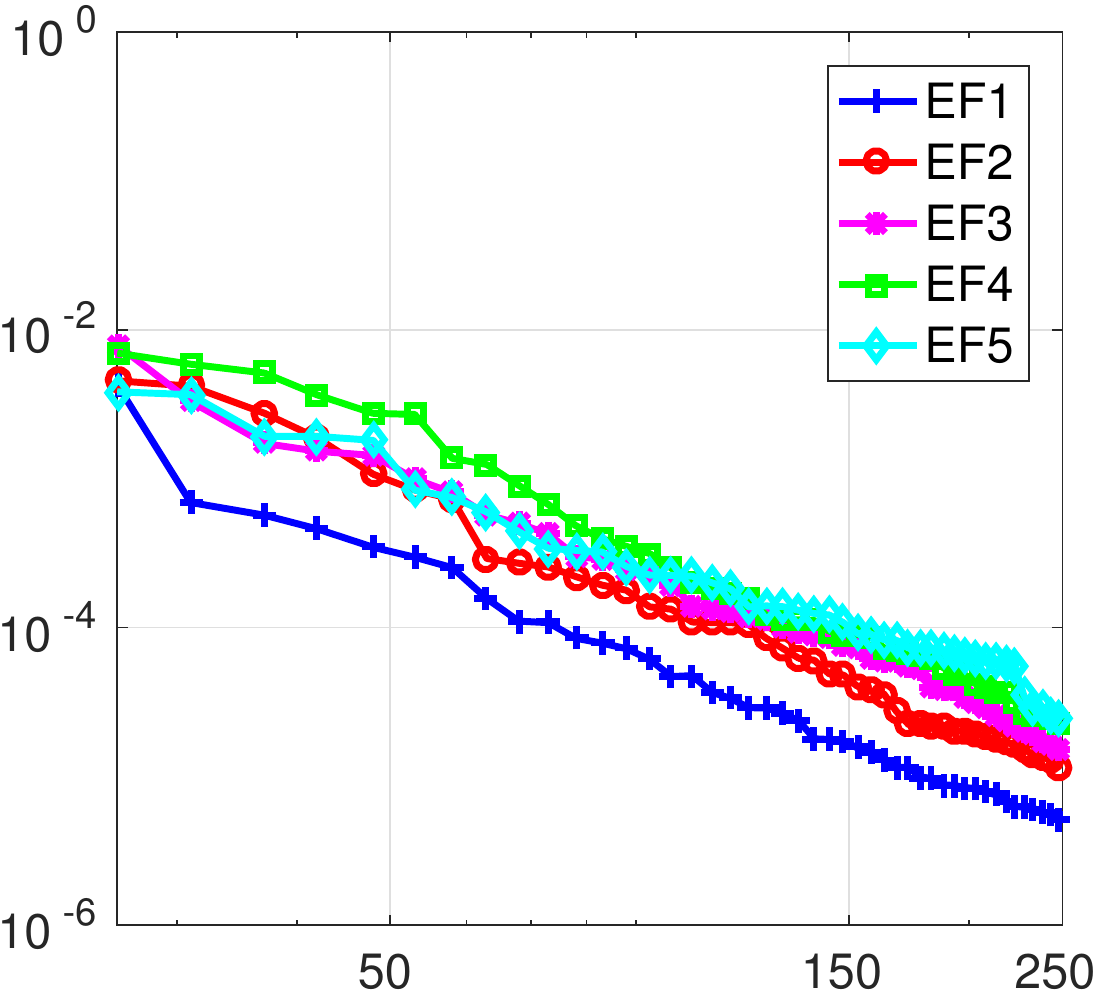}
\includegraphics[width=.32\textwidth]{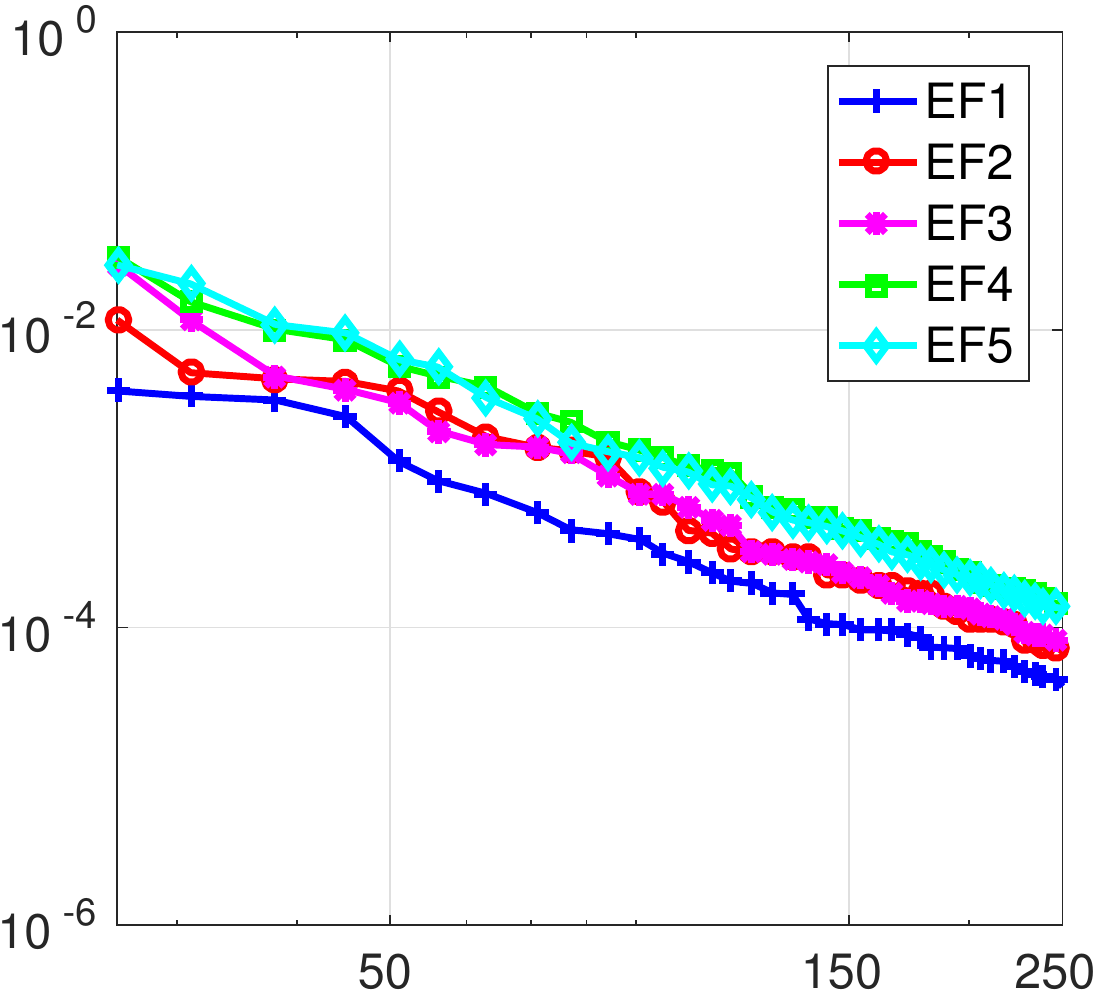}
\includegraphics[width=.32\textwidth]{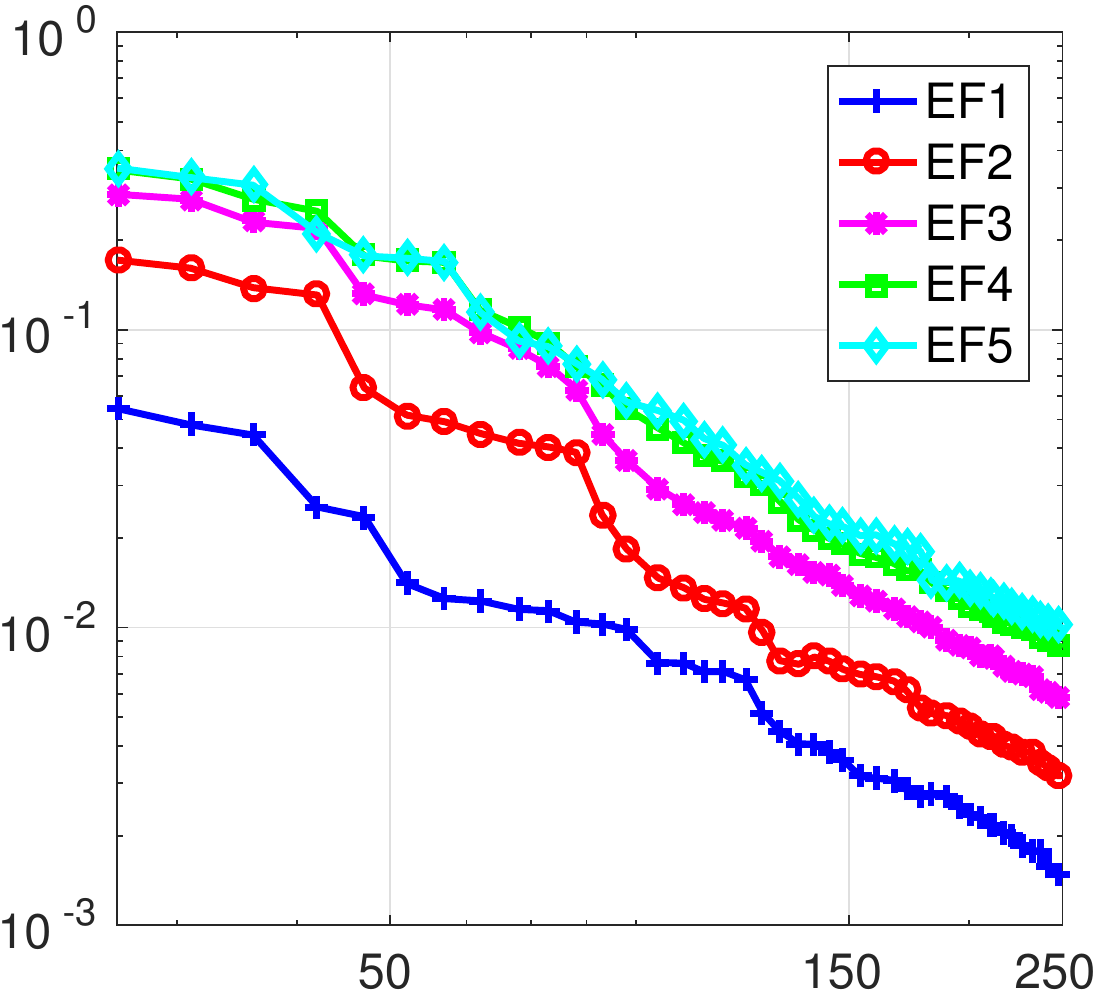}
\caption{Convergence of the relative error of the eigenvalues (top) and eigenfunctions (bottom). Parameter range $\mathcal{P}_1$ with a fixed thickness (left), with varying thickness (middle) and parameter range $\mathcal{P}_2$ with varying thickness (right). }
\label{fig:EW_conv_results}
\end{figure}

The first reduced basis test is the simultaneous approximation of the first five eigenpairs on both parameter-sets $\mathcal{P}_1$ and $\mathcal{P}_2$. We use an initial basis of size $25$ computed by a POD, which is enriched by the greedy algorithm up to a basis size of $250$.  In Figure~\ref{fig:EW_conv_results}, the error decay for the different eigenvalues and eigenfunctions is presented. We observe very good convergence, with a  similar rate in all cases. As expected the magnitude of the error grows with the complexity of the parameter range.

\begin{figure}[htb]
\centering
\includegraphics[width=.45\textwidth]{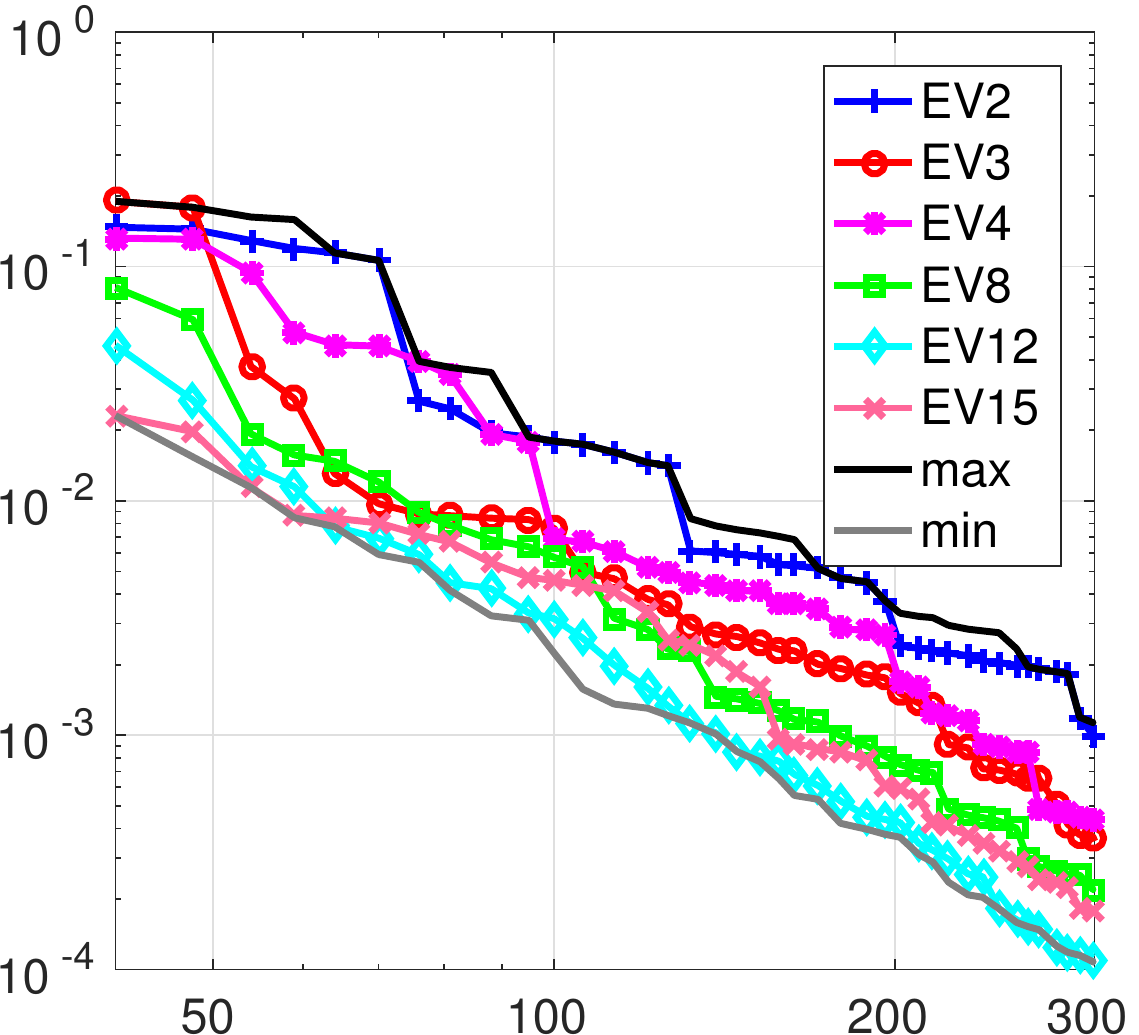} 
\includegraphics[width=.45\textwidth]{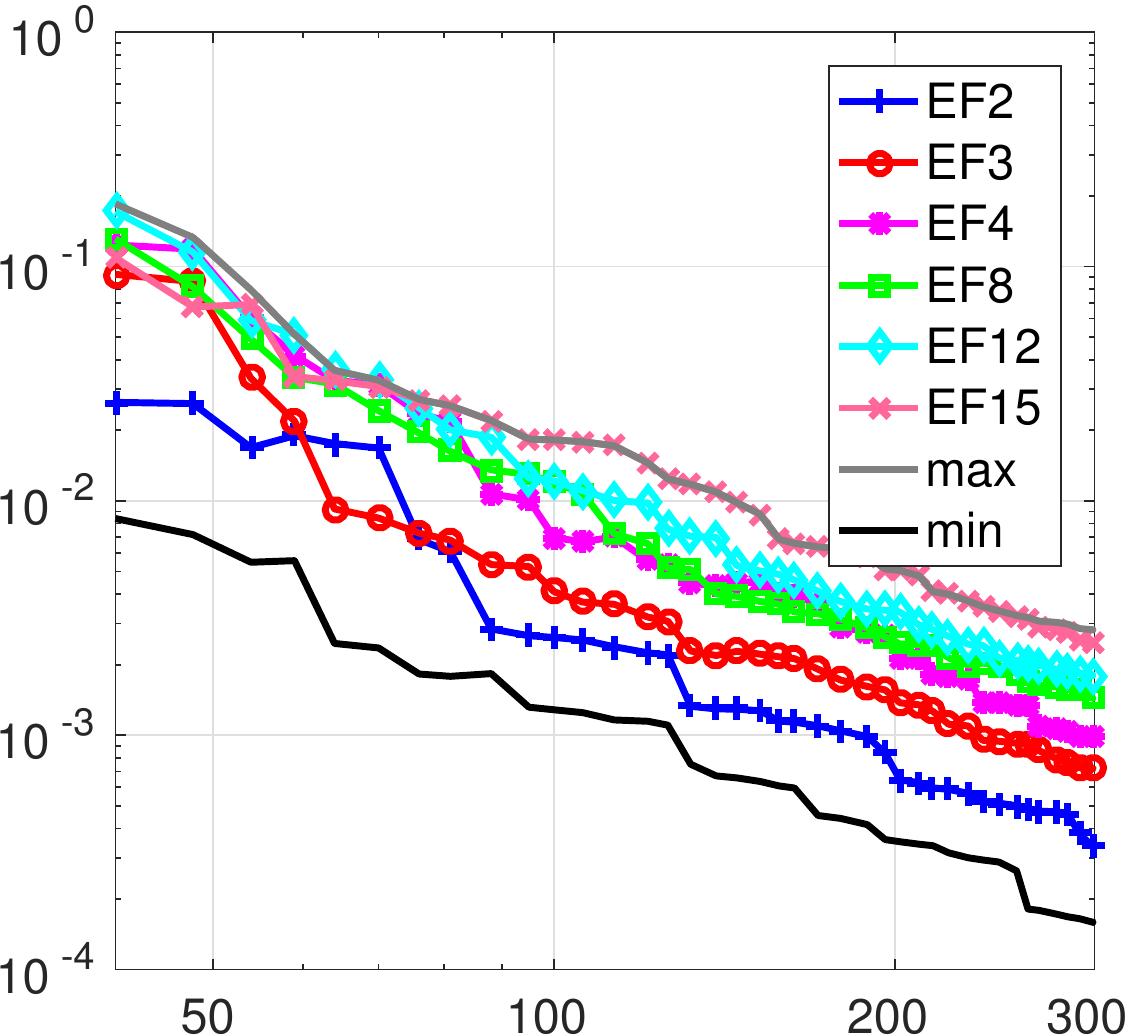}
\caption{Convergence of the relative error of the eigenvalues (left) and eigenfunctions (right). Parameter range $\mathcal{P}_1$ with varying thickness, simultaneous approximating 15 eigenpairs. }
\label{fig:EW_15_conv_results}
\end{figure}

Also an approximation of a larger number of eigenpairs does not pose any unexpected difficulties. Error values for the eigenvalue and eigenfunction are shown in Figure~\ref{fig:EW_15_conv_results} for an approximation of the first 15 eigenpairs in the parameter-set $\mathcal{P}_1$, showing a good convergence behavior. The reduced basis size necessary for a given accuracy increases compared to the previous cases of 5 eigenpairs, due to the higher amount of eigenfunctions which are, for a fixed parameter, orthogonal to each other.

\begin{figure}[htb]
\centering
\begin{minipage}[b]{.48\textwidth}
\includegraphics[width=.85\textwidth]{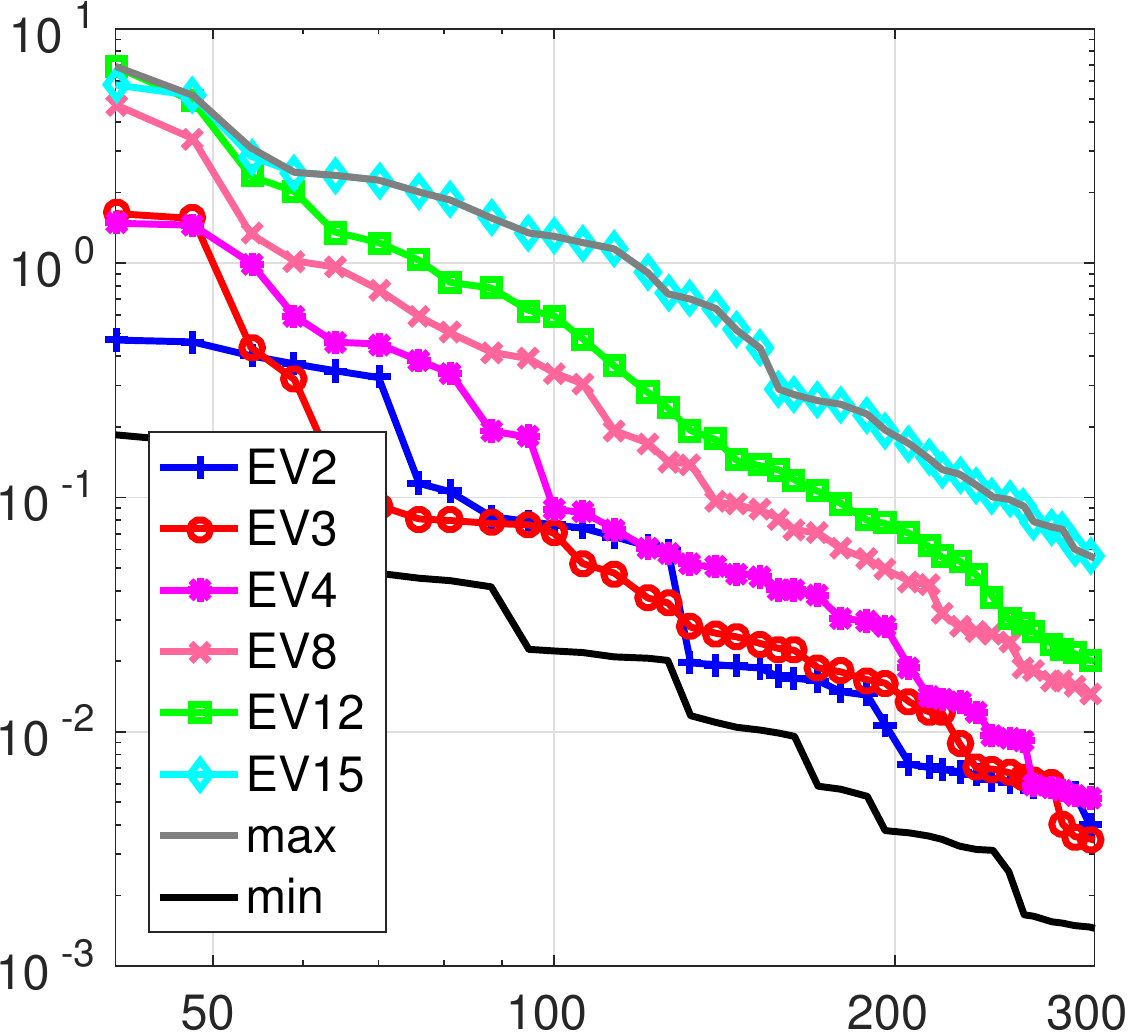}
\caption{Convergence of the absolute error value for the  eigenvalues. Parameter range $\mathcal{P}_1$ with varying thickness, simultaneous approximating 15 eigenpairs.}
\label{fig:EW_15_abs_conv_results}
\vspace*{1em}
\end{minipage}
\hfill
\begin{minipage}[b]{.48\textwidth}
\includegraphics[width=.9\textwidth]{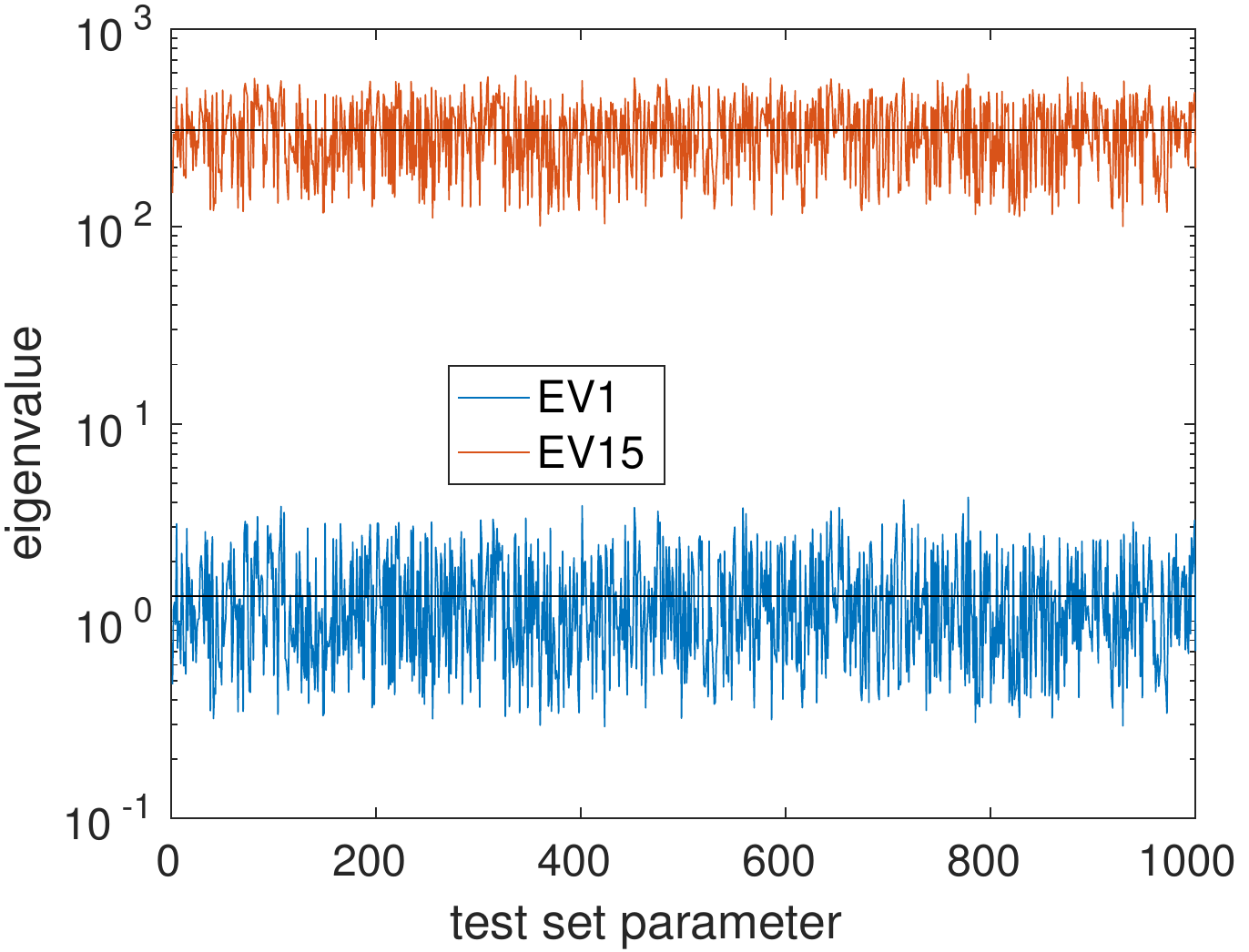}
\caption{ Sampling of the  first and 15th eigenvalue within the parameter set $\mathcal{P}_1$ as used in the test set. Extremal values: $\min \lambda_1 = 0.29$, $\max \lambda_1 = 4.24$, $\min \lambda_{15} = 100.19$, $\max \lambda_{15} = 593.65$. }
\label{fig:ev_values_sampling}
\end{minipage}
\end{figure}

When considering the relative error for the eigenvalues, see Figure~\ref{fig:EW_conv_results} and Figure~\ref{fig:EW_15_conv_results}, we note that for a fixed basis size, the higher eigenvalues have a better relative approximation than the lower ones. In contrast, considering the eigenfunctions,  the error of the ones associated with the lower eigenvalues are smaller compared to the ones associated with the higher eigenvalues. 
In Figure~\ref{fig:EW_15_abs_conv_results}, we see the same effect, when considering absolute error values for the eigenvalue compared to the relative values. The error in the eigenfunctions has the same ordering as the absolute error in the eigenvalue.
For the parameters under consideration, the lower and higher ones of the considered eigenvalues differ by magnitudes as illustrated in Figure~\ref{fig:ev_values_sampling}. Hence when computing the relative error from the absolute ones, the error values of the high eigenvalues are divided by a large number and become small compared to the lower eigenvalues.

\newpage

\section{Conclusion}
\label{sec:Conclusion}
We have shown an eigenvalue reduced basis approximation of a violin bridge as an interesting application in vibroacoustic. 
The model reduction yields a good  approximation  quality in all cases of consideration, with a significant complexity reduction. The detailed system, a saddle point problem of $47,985$ degrees of freedom, was reduced to a positive-definite system of less than $300$ degrees of freedom. Furthermore, the thickness of the violin bridge, included as a geometry parameter, was shown to have a significant influence on the eigenvalues and eigenfunctions, without posing further difficulties to the reduced basis approximation. Altogether, it was shown that reduced basis methods are suitable for vibroacoustical mortar settings even with a parameter-dependent geometry.
\newpage

\bibliographystyle{siam}
\bibliography{./lit}
\end{document}